\documentclass[preprint,12pt]{elsarticle}

\usepackage{amssymb}
\usepackage{amsmath}
\usepackage{accents}
\usepackage{graphicx}
\usepackage[colorinlistoftodos]{todonotes}
\usepackage[colorlinks=true, allcolors=black]{hyperref}
\usepackage{lipsum} % for dummy text only
\usepackage{floatrow} % for figures next to each other
\usepackage[hypcap=false]
{caption} % for caption of tables
\usepackage{subcaption}

\usepackage{listings}
\usepackage{url}
\usepackage{courier}
\usepackage{wrapfig}
\usepackage{color,soul}
\definecolor{codegreen}{rgb}{0,0.6,0}
\usepackage{chngpage}
\usepackage{xcolor}
\usepackage{bm}
\usepackage{amsmath,amssymb}
\usepackage{verbatim}
\newlength{\mycolwidth}
\usepackage{array}
\newcolumntype{Z}{>{$}p{\mycolwidth}<{$}}
\usepackage{algorithmicx}
\usepackage{algorithm} % for algorithms
\usepackage{refcount}
\usepackage{algpseudocode} % for algorithms 2
\usepackage{placeins} % for \Floatbarrier
\MakeRobust{\Call}

\usepackage{amsfonts}
\usepackage{graphicx}
\usepackage{epstopdf}
\usepackage{amsmath} % useful math package

\usepackage{lscape}
\usepackage{booktabs}

\usepackage{tikz,lipsum,lmodern}
\usepackage[most]{tcolorbox}

\usepackage{pifont}% http://ctan.org/pkg/pifont
\newcommand{\cmark}{\textcolor{blue}{\ding{51}}}%
\newcommand{\xmark}{\textcolor{gray}{\ding{55}}}

\DeclareMathAccent{\svec}{\mathord}{letters}{126}
\newcommand\acclrvec[1]{\accentset{\,\leftrightarrow}{#1}}
\newcommand\stvec[1]{\mathbf #1}
		
\newcommand\ssvec[1]{\acclrvec{\stvec{#1}}}	

\newcommand\cssvec[1]{\acclrvec{\tilde{\stvec{#1}}}}

\usepackage{scalerel,stackengine}
\stackMath

\newcounter{bla}

\journal{Computer Physics Communications}

\begin{document}

\begin{frontmatter}



\title{{\fontfamily{qcr}\selectfont  HORSES3D}: a high-order discontinuous Galerkin solver for flow simulations and multi-physics applications}


\author[1,2]{E. Ferrer}
\author[1,2]{G. Rubio \corref{cor}}
%\ead{g.rubio@upm.es}
\author[1]{G. Ntoukas}
\author[1]{W. Laskowski}
\author[1]{O.A. Mariño}
\author[1]{S. Colombo}
\author[1]{A. Mateo-Gab\'in}
\author[1]{F. Manrique de Lara}
\author[1]{D. Huergo}
\author[5]{J. Manzanero}
\author[3]{A.M. Rueda-Ramírez}
\author[4]{D.A. Kopriva}
\author[1,2]{E. Valero}

%\author[1,2]{E \snm{Ferrer} et al.}
\cortext[author] {Corresponding author.\\\textit{E-mail address:} g.rubio@upm.es}
\address[1]{ETSIAE-UPM-School of Aeronautics, Universidad Politécnica de Madrid, Madrid-Spain}
\address[2]{Center for Computational Simulation, Universidad Politécnica de Madrid, Madrid-Spain}
\address[3]{University of Cologne, Division of Mathematics, Köln-Germany}
\address[4]{Department of Mathematics, Florida State University and Computational Science Research Center, San Diego State University, USA}
\address[5]{Airbus Defence and Space, Madrid-Spain}

\begin{abstract}
We present the latest developments of our \textbf{H}igh-\textbf{O}rder \textbf{S}pectral \textbf{E}lement \textbf{S}olver ({\fontfamily{qcr}\selectfont  HORSES3D}), an open source high-order discontinuous Galerkin framework, capable of solving a variety of flow applications, including compressible flows (with or without shocks), incompressible flows, various RANS and LES turbulence models, particle dynamics, multiphase flows, and aeroacoustics.
We provide an overview of the high-order spatial discretisation (including energy/entropy stable schemes) and anisotropic p-adaptation capabilities. The solver is parallelised using MPI and OpenMP showing good scalability for up to 1000 processors. Temporal discretisations include explicit, implicit, multigrid, and dual time-stepping schemes with efficient preconditioners. Additionally, we facilitate meshing and simulating complex geometries through a mesh-free immersed boundary technique. 
We detail the available documentation and the test cases included in the GitHub repository.
%the different models implemented for turbulent flows (ILES, Smagorinsky, SVV-Smagorinsky, WALE) and multiphase flows (NS+Cahn-Hilliard). Additionally, we detail the capabilities to perform local p-adaptation, and the preconditioners (including multigrid) to advance the solution in time, efficiently.
\end{abstract}




\end{frontmatter}

% All CPiP articles must contain the following
% PROGRAM SUMMARY.

{\bf PROGRAM SUMMARY}

\begin{small}
\noindent
{\em Program Title: } {\fontfamily{qcr}\selectfont  HORSES3D}                                         \\
{\em CPC Library link to program files:} (to be added by Technical Editor) \\
{\em Developer's repository link:} private GitHub repository. \\
{\em Code Ocean capsule:} (to be added by Technical Editor)\\
{\em Licensing provisions:} MIT  \\
{\em Programming language:} Fortran 2008                                  \\
{\em External routines/libraries:} METIS, MPI, HDF5, MKL, PETSc (all are optional).                                  \\
{\em Nature of problem:} {\fontfamily{qcr}\selectfont  HORSES3D} is a high-order discontinuous Galerkin framework, capable of solving a variety of flow applications, including compressible flows (with or without shocks), incompressible flows, various RANS and LES turbulence models, particle dynamics, multiphase flows, and aeroacoustics.\\
  %Describe the nature of the problem here. \\
{\em Solution method:} high-order discontinuous Galerkin Spectral Element (DGSEM) and explicit/implicit time-steppers.\\
  %Describe the method solution here.
{\em Additional comments including restrictions and unusual features:} {\fontfamily{qcr}\selectfont  HORSES3D} is a multiphysics environment where the compressible Navier-Stokes equations, the incompressible Navier–Stokes equations, the Cahn–Hilliard equation and entropy–stable variants are solved. Arbitrary high–order, p–anisotropic discretisations are used, including static and dynamic p–adaptation methods (feature-based and truncation error-based). Explicit and implicit time-steppers for steady and time-marching solutions are available, including efficient multigrid and preconditioners. Numerical and analytical Jacobian computations with a coloring algorithm have been implemented.
Multiphase flows are solved using a diffuse interface model: Navier–Stokes/Cahn–Hilliard.
Turbulent models implemented include RANS: Spalart-Allmaras and LES: Smagorinsky, Wale, Vreman; including wall models.
Immersed boundary methods can be used, to avoid creating body fitted meshes.
Acoustic propagation can be computed using Ffowcs-Williams and Hawkings models.
{\fontfamily{qcr}\selectfont  HORSES3D} supports curvilinear, hexahedral, conforming meshes in GMSH, HDF5 and SpecMesh/HOHQMesh format.
A hybrid CPU-based parallelisation strategy (shared and distributed memory) with OpenMP and MPI is followed.
  %Provide any additional comments here.

\end{small}

 \tableofcontents

% \linenumbers

%\section{Program Summary}
%\textbf{Program Title}: {\fontfamily{qcr}\selectfont  HORSES3D}

%\textbf{Program Files doi}: 

%\textbf{Code Ocean Capsule}: 

%\textbf{Licensing provisions}: MIT

%\textbf{Programming language}: Fortran 2008

%\textbf{External routines/libraries}: METIS, MPI, HDF5, MKL, PETSc (all are optional).

%\textbf{Nature of problem}: {\fontfamily{qcr}\selectfont  HORSES3D} is a high-order discontinuous Galerkin framework, capable of solving a variety of flow applications, including compressible flows (with or without shocks), incompressible flows, various RANS and LES turbulence models, particle dynamics, multiphase flows, and aeroacoustics.

%\textbf{Solution method}: high-order discontinuous Galerkin Spectral Element (DGSEM) and explicit/implicit time-steppers.

\section{Introduction}
\label{sec:intro}
To simulate and predict complex flows, the non-linear Navier-Stokes (NS) equations that govern fluid flow need to be solved numerically. Traditional simulation solvers that use low-order methods (order $\leq$ 2, for example, finite differences or finite volumes) have been widely adopted by industry and academia to simulate fluid flows, but have limitations when high accuracy is required \cite{WorkshopDG}. Low-order numerical techniques (most commercial codes) suffer from non-negligible numerical errors (e.g., dissipative and dispersive errors) that can increase unphysical dissipation and dispersion of flow structures and provide unrealistic results. 

An alternative, to low order solvers, is to employ high-order discretisations based on polynomial approximations of the solution. These methods are characterised by low numerical errors, and their ability to use mesh refinement through an increased number of mesh nodes (h-refinement) and/or polynomial enrichment (p-refinement) to achieve highly accurate solutions. Such high-order polynomial methods produce an exponential decay of the error for smooth solutions instead of the algebraic decay characteristic of low-order techniques. 

Typical aircraft aerodynamic simulations require an order of 50 to 100M degrees of freedom ($dof$), if low-order methods are used. Using high-order methods, the overall number of $dof$ can be reduced for a desired level of accuracy. 
Given a numerical scheme, the numerical error $e$ is proportional to $(dof)^P$, where $P$ is the order, see \cite{wang2016perspective}. To achieve similar errors (i.e. $e_{HO} \approx e_{LO}$) between high- ($HO$) and low-order ($LO$) methods, the high-order mesh may be coarsened following: $(dof)_{HO}=\exp(P_{LO}/P_{HO}.\log[(dof)_{LO}])$. Selecting $P_{LO}=2$ and $P_{HO}=5$, one finds for $(dof)_{LO}=100M$ a corresponding mesh size $(dof)_{HO}=1.6k$; that is, a factor of 60000 decrease in $dof$ for the same accuracy. Of course, this estimate is optimistic, since in reality small cells are required near walls to resolve boundary layers. Nevertheless, this estimate illustrates how high-order methods can enable simulations with fewer degrees of freedom (while maintaining accuracy). Additional accuracy can be achieved with a limited increase in cost, when using local anisotropic p-refinement (i.e. only particular mesh elements and directions use high-order polynomials), \cite{kompenhans2016adaptation,kompenhans2016comparisons,rueda2019truncation,rueda2019p}.

High-order numerical methods of spectral type (order $>$ 2) were first developed in the 1970s (e.g., \cite{gottlieb1977numerical}) and have typically been perceived as accurate but expensive methods, with a lack of flexibility (for complex geometries) and robustness (tendency to crash) when applied to complex turbulent flows. Although widely adopted in academia and particularly for low Reynolds numbers during the 1980s and 90s, it was not until the 2000s that these became a realistic alternative to well-established low-order techniques (finite differences or finite volumes). One of the main reasons for the increasing acceptance has been the development of high-order discontinuous Galerkin (DG) variants that include interelement fluxes that enhance robustness (when compared to classic high-order continuous methods). DG has seen a remarkable recent increase in popularity in the last 10 years. DG methods were first developed by Reed $\&$ Hill in the framework of neutron transport equations \cite{reed1973triangular}, but remained unexploited until the late 1990s, when the method was generalised to convection-diffusion problems that are used to simulate fluid flows: see, e.g., \cite{bassi1997high,bassi1997high2}. Today, DG methods have been used to solve a wide variety of flows, including incompressible, compressible, or multiphase flows, 
%see (Ferrer et al. 2014), (Ferrer $\&$ Willden, 2015), (Kompenhans et al. 2016), (Kompenhans et al. 2016b), (Manzanero et al. 2017), (Ferrer, 2017), (Manzanero et al. 2018), 
but still need developments to become accepted within the industry (e.g., increased speed and robustness for turbulent flows). 

DG is paving its way towards industry through research centres. For example, Airbus, Dyson, McLaren F1 or Siemens-Gamesa have collaborated with our group in various European projects to evaluate high-order methods. In these projects, aeronautical European research centres (e.g., ONERA or DLR) were also involved.  DG methods are relatively mature for aeronautical applications, but require further developments to be widely adopted by industry (not necessarily aeronautical).

In this work, we summarise the main features included in our high-order DG solver {\fontfamily{qcr}\selectfont  HORSES3D}:

\begin{tcolorbox}[colback=black!5!white,colframe=gray!75!black]
\begin{itemize}

\item     A multiphysics environment where the compressible Navier-Stokes equations, the incompressible Navier–Stokes equations, the Cahn–Hilliard equation and entropy–stable variants are solved.
\item     Arbitrary high–order, p–anisotropic discretisations, including static and dynamic p–adaptation methods using feature-based and truncation error estimates.
\item     Explicit and implicit time-steppers for steady and time-marching solutions, including efficient multigrid and preconditioners.
%\item     Numerical and analytical Jacobian computations with a coloring algorithm.
\item Multiphase flows using a diffuse interface model: Navier–Stokes/Cahn– Hilliard.
\item     Turbulent models including RANS: Spalart-Allmaras and LES: Smagorinsky, Wale, Vreman; including wall models.
\item     Immersed boundary methods to avoid creating body fitted meshes.
\item     Acoustic propagation using Ffowcs-Williams and Hawkings.
\item     Support for curvilinear, hexahedral, conforming meshes in GMSH, HDF5 and SpecMesh/HOHQMesh format.
\item     Hybrid CPU-based parallelisation (shared and distributed memory) with OpenMP and MPI.
\end{itemize}
\end{tcolorbox}

The work is structured as follows. In Sec.~\ref{sec:HOsolvers}, we contextualise {\fontfamily{qcr}\selectfont  HORSES3D} with a review of the main capabilities of available high-order academic solvers. The rest of the paper can be split in two parts. 
In the first part: Numerical Method 
(Sec.~\ref{sec:themethods}), we describe the original features of the solver that was first conceived to solve the compressible Navier-Stokes equations. This part includes meshing (Sec.~\ref{sec:mesh}), the compressible Navier-Stokes solver  (Sec.~\ref{sec:cNS}), the discontinuous Galerkin spatial discretisation (Sec.~\ref{sec:dg}), the possibility to include alternate physics (Sec.~\ref{sec:altenPhy}), the available energy/entropy stable formulations to enhance stability (Sec.~\ref{sec:stable}), local spatial p-adaption (Sec.~\ref{sec:space_adapt}), the time advancement methods available (Sec.~\ref{sec:time}) and parallelisation benchmarks (Sec.~\ref{sec:parallel_cost}).

In the second part: Multiphysics (Sec.~\ref{sec:multiphysics}), we detail the extended capabilities of {\fontfamily{qcr}\selectfont  HORSES3D} , including turbulent modelling (Sec.~\ref{sec:turb}), incompressible flows (Sec.~\ref{sec:iNS}), multiphase flows (Sec.~\ref{sec:Mul}), shock capturing (Sec.~\ref{sec:Shock}), particle dynamics (Sec.~\ref{sec:particles}), aeroacoustics (Sec.~\ref{sec:acoustic}) and immersed boundaries (Sec.~\ref{sec:IB}).

\section{Review of high-order solvers}
\label{sec:HOsolvers}
A range of high-order academic solvers exist and are briefly summarised in Table \ref{tab:my-table}, where we include their main capabilities. We provide references for each solver and outline here the capabilities of interest when compared to {\fontfamily{qcr}\selectfont  HORSES3D}. This list is a non-exhaustive summary where we have not included commercial software. The table reflects that {\fontfamily{qcr}\selectfont  HORSES3D} provides a unique set of capabilities to simulate complex flows.
%\cite{CANTWELL2015205,MOXEY2020107110}
% \usepackage{booktabs}
% \usepackage{graphicx}
% \usepackage{booktabs}
% \usepackage{graphicx}
% \usepackage{booktabs}
% \usepackage{graphicx}
% \usepackage{lscape}
% \usepackage{booktabs}
% \usepackage{lscape}
% \usepackage{booktabs}
% \usepackage{lscape}
%\begin{landscape}
\begin{table}[]
\centering
\resizebox{\textwidth}{!}{%
\begin{tabular}{@{}ccccccccc@{}}
\toprule
%\footnotesize
\textbf{} &
  \textbf{{\fontfamily{qcr}\selectfont  HORSES3D}} &
  \textbf{Nek5000} &
  \textbf{Nektar++} &
  \textbf{Semtex} &
  \textbf{deal.II} &
  \textbf{Flexi/Fluxo} &
  \textbf{Trixi} &
  \textbf{PyFR} \\ 
  %\hline
  \midrule
\textbf{References} &
   &
   \cite{Nek5000}&
  \cite{CANTWELL2015205,MOXEY2020107110} &
  \cite{BLACKBURN2019106804}   &
  \cite{10.1145/1268776.1268779} &
  \cite{GASSNER201639,HINDENLANG201286} &
  \cite{ranocha2022adaptive,schlottkelakemper2021purely}&
   \cite{WITHERDEN20143028}\\
\textbf{Language} &
  Fortran &
  Fortran &
  C++ &
  C++ &
  C++ &
  Fortran &
  Julia &
  Python \\
\textbf{Spatial discretisation} &
  DG &
  CG &
  CG/DG &
  CG &
  CG/DG &
  DG &
  DG &
  FR \\
\textbf{Incompressible} &
  \cmark &
  \cmark &
   \xmark &
  \cmark &
  \cmark &
   \xmark&
  \xmark &
  \cmark \\
\textbf{Compressible} &
  \cmark &
   \xmark&
  \cmark &
  \xmark &
  \cmark &
  \cmark &
  \cmark &
  \cmark \\
%\textbf{\begin{tabular}[c]{@{}c@{}}Energy/Entropy \\ stable\end{tabular}} &
\textbf{Energy/Entropy stable} &
  \cmark &
  \xmark &
  \xmark &
  \xmark &
  \xmark &
  \cmark &
  \cmark &
  \xmark \\
\textbf{p-adaptation} &
  \cmark &
  \cmark &
  \cmark &
  \xmark &
  \cmark &
  \cmark &
  \cmark &
  \cmark \\
%\textbf{\begin{tabular}[c]{@{}c@{}}Steady/Unsteady\\ in time \end{tabular}} &
\textbf{Steady/Unsteady} &
  S/U &
  U &
  U &
  U &
  S/U &
  U &
  U &
  S/U \\
\textbf{LES} &
  \cmark &
  \cmark &
  \cmark &
  \cmark &
  \cmark &
  \cmark &
  \cmark &
  \cmark \\
\textbf{RANS} &
  \cmark &
   \xmark&
   \xmark&
  \xmark &
  \xmark &
  \xmark &
  \xmark &
  \xmark \\
\textbf{Multiphase} &
  \cmark &
  \cmark &
  \xmark &
  \xmark &
  \cmark &
  \xmark &
  \xmark &
  \cmark \\
%\textbf{\begin{tabular}[c]{@{}c@{}}Shock \\ Capturing\end{tabular}} &
\textbf{Shocks} &
  \cmark &
  \xmark &
  \cmark &
  \xmark &
  \cmark &
  \cmark &
  \cmark &
  \cmark \\
%\textbf{\begin{tabular}[c]{@{}c@{}}Particle \\ Dynamics\end{tabular}} &
\textbf{Particles} &
  \cmark &
  \cmark &
  \cmark &
  \cmark &
  \cmark &
  \xmark &
  \xmark &
  \xmark \\
\textbf{Aeroacoustics} &
  \cmark &
  \xmark &
  \cmark &
  \xmark &
  \cmark &
  \xmark &
  \cmark &
  \xmark \\
%\textbf{\begin{tabular}[c]{@{}c@{}}Mesh Free \\ Immersed Boundaries\end{tabular}} &
\textbf{Immersed Boundaries} &
  \cmark &
  \xmark &
  \cmark &
  \xmark &
  \xmark &
  \xmark &
  \xmark &
   \xmark \\
\multicolumn{1}{l}{} &
  \multicolumn{1}{l}{} &
  \multicolumn{1}{l}{} &
  \multicolumn{1}{l}{} &
  \multicolumn{1}{l}{} &
  \multicolumn{1}{l}{} &
  \multicolumn{1}{l}{} &
  \multicolumn{1}{l}{} &
  \multicolumn{1}{l}{} \\ \bottomrule
\end{tabular}
}
\caption{Summary of high-order solvers and capabilities.}
\label{tab:my-table}
\end{table}
%\end{landscape}

\section{Numerical methods}
\label{sec:themethods}

{\fontfamily{qcr}\selectfont  HORSES3D} is a Fortran 2003 object-oriented solver, originally developed to solve the compressible Navier-Stokes equations using DG discretisations and explicit time marching. 
In this section, we detail the compressible solver together with the spatial- and temporal-time advancement schemes.

\subsection{Meshing and Post-processing}
\label{sec:mesh}
The solver supports curvilinear, hexahedral and conforming meshes in GMSH \cite{geu09}, HDF5\cite{folk1999hdf5} and SpecMesh/HOHQMesh \cite{HOHQMesh} format. Linear meshes generated in Gambit, ICEM, or CGNS can be curved and converted to HDF5 using HOPR \cite{HOPR}.  If the user prefers to avoid generating body fitted meshes, we provide a mesh-free immersed boundary capability, which is detailed in (Sec.~\ref{sec:IB}).

Tecplot and Paraview \cite{ahrens2005paraview,ayachit2015paraview} can be used to post-process the results. We have developed a post-processing utility \textit{horses2plt} that allows interpolating to homogeneous nodal points (from Gauss-Legendre or Gauss-Lobatto) and modifying the polynomial order for visualisation purposes. The tool also allows to extract derived variables, such as vorticity and Q-criterion (if gradients are enabled during the simulation) or Reynolds Stresses (if the statistics are enabled during the simulation).

\subsection{Compressible Navier-Stokes}
\label{sec:cNS}
 
The 3D Navier-Stokes equations can be compactly written as
% To do so, we consider the DG approximation of the compressible Navier-Stokes equations:
%
\begin{equation}
\boldsymbol{u}_t+\nabla\cdot\ssvec{{F}}_e = \nabla\cdot\ssvec{F}_v,
\label{eq:compressibleNScompact}
\end{equation}
where $\boldsymbol{u}$ is the state vector of conservative variables $\boldsymbol{u} = [ \rho , \rho v_1 , \rho v_2 , \rho v_3 , \rho e]^T$, $\ssvec{F}_e$ are the inviscid, or Euler fluxes,
\begin{equation}
\ssvec{F}_e = \left[\begin{array}{ccc} \rho v_1 & \rho v_2 & \rho u_3 \\
                                                                                \rho v_1^2 + p & \rho v_1v_2 & \rho v_1v_3 \\
                                                                                	\rho v_1v_2 & \rho v_2^2 + p & \rho v_2v_3 \\
                                                                                	\rho v_1v_3 & \rho v_2v_3 & \rho v_3^2 + p \\
                                                                                	\rho v_1 H & \rho v_2 H & \rho v_3 H
\end{array}\right],
\end{equation}
where $\rho$, $e$, $H=E+p/\rho$, and $p$ are the density, total energy, total enthalpy and pressure, respectively, and $\vec{v}=[v_1,v_2,v_3]^T$ is the velocity. Additionally, $\ssvec{F}_v$ defines the viscous fluxes,
\begin{equation}
\ssvec{F}_v= \left[\begin{array}{ccc}0 & 0 & 0\\
 \tau_{xx} & \tau_{xy} & \tau_{xz} \\
 \tau_{yx} & \tau_{yy} & \tau_{yz} \\
 \tau_{zx} & \tau_{zy} & \tau_{zz} \\
 \sum_{j=1}^3 v_j\tau_{1j} + \kappa T_x& \sum_{j=1}^3 v_j\tau_{2j} + \kappa T_y& \sum_{j=1}^3 v_j\tau_{3j} + \kappa T_z
\end{array}\right],
\label{eq:viscousfluxes}
\end{equation}
where $\kappa$ is the thermal conductivity, $T_x, T_y$ and $T_z$ denote the temperature gradients and the stress tensor $\boldsymbol{\tau}$ is defined as $\boldsymbol{\tau} = \mu(\nabla \vec{v} + (\nabla \vec{v})^T) - 2/3\mu \boldsymbol{I}\nabla\cdot\vec{v}$, with $\mu$ the dynamic viscosity and $\boldsymbol{I}$ the three-dimensional identity matrix.

\subsection{Spatial discretisation: discontinuous Galerkin}
\label{sec:dg}

{\fontfamily{qcr}\selectfont  HORSES3D} discretises the Navier-Stokes equations using the Discontinuous Galerkin Spectral Element Method (DGSEM), which is a particularly efficient nodal version of DG schemes \cite{2009:Kopriva}. For simplicity, here we only introduce the fundamental concepts of DG discretisations. More details can be found in \cite{MANZANERO2020109241,FerrerJCP}.

% The 3D Navier-Stokes equations can be  written as:
% % To do so, we consider the DG approximation of the compressible Navier-Stokes equations:
% %
% \begin{equation}\label{eq::NS1}
% \vec{u}_t   +\nabla\cdot\vec{\boldsymbol{F}}_e= \nabla\cdot\vec{\boldsymbol{F}}_v,
% \end{equation}
% %
% where $\vec{u}$ is the vector of conservative variables $\vec{\boldsymbol{u}} = ( \rho , \rho v_1 , \rho v_2 , \rho v_3 , \rho e)^T$ in compressible solvers. For incompressible solvers $\vec{\boldsymbol{u}} = (  v_1 ,  v_2 ,  v_3 )^T$ and Eq.\eqref{eq::NS1} is complemented with $\nabla\cdot \boldsymbol{u}$. Details on the definition of inviscid and viscous solvers can be found in  \cite{manzanero2018,FerrerJCP}.
%Details on the specific formulations retained for inviscid and viscous fluxes $\vec{\boldsymbol{F}}_e$ and $\vec{\boldsymbol{F}}_v$ can be found in \ref{appendix::A} of this text.
%
The physical domain is tessellated with non-overlapping curvilinear hexahedral elements, $e$, which are geometrically transformed to a reference element, $el$. This transformation is performed using a polynomial transfinite mapping that relates the physical coordinates $\vec{x}$ and the local reference coordinates $\vec{\xi}$. The transformation is applied to \eqref{eq:compressibleNScompact} resulting in the following:

\begin{equation}
J \boldsymbol{u}_t  + \nabla_\xi\cdot\cssvec{F}_e = \nabla_\xi\cdot\cssvec{F}_v,
\label{eq:compressibleNScompact_transformed}
\end{equation}

where $J$ is the Jacobian of the transfinite mapping, $\nabla_\xi$ is the differential operator in the reference space and $\cssvec{F}$ are the contravariant fluxes \cite{2009:Kopriva}. 

To derive DG schemes, we multiply \eqref{eq:compressibleNScompact_transformed} by a locally smooth test function $\phi_j$, for $0\leq j\leq P$, where $P$ is the polynomial degree, and integrate over an element $el$ to obtain the weak form
\begin{equation}\label{eq::NS2}
\int_{el}J \boldsymbol{u}_t\phi_j+\int_{el} \nabla_\xi\cdot\cssvec{F}_e\phi_j  =\int_{el} \nabla_\xi\cdot\cssvec{F}_v\phi_j.
\end{equation}
We can now integrate by parts the term with the inviscid fluxes, $\cssvec{F}_e$, to obtain a local weak form of the equations (one per mesh element) with the boundary fluxes separated from the interior
\begin{equation}\label{eq::NS3}
\int_{el}J \boldsymbol{u}_t\phi_j +  \int_{\partial el} \cssvec{F}_e\cdot\hat{\mathbf{n}}\phi_j-\int_{el} \cssvec{F}_e\cdot\nabla_\xi\phi_j
=\int_{el} \nabla_\xi\cdot\cssvec{F}_v\phi_j,
\end{equation}
where $\hat{n}$ is the unit outward vector of each face of the reference element ${\partial el}$. %Summing over all element in the mesh,
We replace discontinuous fluxes at inter--element faces by a numerical inviscid flux, $\cssvec{F}_{e}^{\star}$, to couple the elements, %(that includes all elements $el$ in the mesh $\Omega$ and element faces $\partial el$ in the set of faces in the mesh $\Gamma$).
\begin{equation}\label{eq::NS4}
\int_{el}J \boldsymbol{u}_t\cdot\phi_j + \int_{\partial el} {\cssvec{F}_{e}^{\star}}\cdot\hat{\mathbf{n}}\phi_j-\int_{el} \cssvec{F}_e\cdot\nabla_\xi\phi_j
=\int_{el} \nabla_\xi\cdot\cssvec{F}_v\phi_j.
\end{equation}
This set of equations for each element is coupled through the inviscid fluxes $\cssvec{F}_{e}^{\star}$ and governs the numerical characteristics, see for example the classic book by Toro \cite{toro2009riemann}. Note that one can proceed similarly and integrate by parts the viscous terms (see, for example, \cite{Unified,ferrer2016a,FERRER2011224,FERRER20127037}). The viscous terms require further manipulations to obtain usable discretisations (Bassi Rebay 1 and 2 or Interior Penalty) and, for the sake of simplicity, here we retain the simple volume form:
\begin{equation}\label{eq::NS5}
\int_{el}J \boldsymbol{u}_t\cdot\phi_j + \int_{\partial el} \underbrace{{\cssvec{F}_{e}^{\star}}\cdot\hat{\mathbf{n}}}_\text{Convective fluxes}\phi_j-\int_{el} \cssvec{F}_e\cdot\nabla_\xi\phi_j
=\int_{el} (\underbrace{\nabla_\xi\cdot\cssvec{F}_v}_\text{Viscous term})\cdot\phi_j.
\end{equation}

%\begin{equation}\label{eq::NS5}
%\int_{el}\vec{u}_t\cdot\phi_j + \int_{\partial el} \underbrace{{\vec{\boldsymbol{F}}_e^*}\cdot\mathbf{n}}_\text{Convective fluxes}\phi_j-\int_{el} \vec{\boldsymbol{F}}_e\cdot\nabla\phi_j
%=\underbrace{\int_{\partial el} \vec{\boldsymbol{F}}_v\cdot\mathbf{n}\phi_j
%-\int_{el} {\vec{\boldsymbol{F}}_v}\cdot\nabla\phi_j}_\text{Viscous terms}.
%\end{equation}

The final step, to obtain a usable numerical scheme, is to approximate the numerical solution and fluxes by polynomials (of order $P$) and to use Gaussian quadrature rules to numerically approximate volume and surface integrals. In {\fontfamily{qcr}\selectfont  HORSES3D} we allow for Gauss-Legendre and Gauss--Lobatto quadrature points. The former are more accurate but the latter allow for the derivation of energy/entropy stable formulations, see details in (Sec.~\ref{sec:stable}) %The non-linear inviscid and viscous terms that can be discretised to control dissipation in the numerical scheme have been underlined.

\begin{tcolorbox}[colback=black!5!white,colframe=gray!75!black]
  In {\fontfamily{qcr}\selectfont  HORSES3D} we have implemented the following convective fluxes: 
Central, Lax-Friedrichs, Rusanov, Roe, Low dissipation Roe, Roe-Pike, matrix dissipation; 
and the viscous fluxes: Bassi-Rebay 1 and 2 and Interior Penalty (symmetric version). 
\end{tcolorbox}
Riemann solvers are the classic option to include numerical dissipation in DG schemes \cite{Gassner_ILES,WorkshopDG}, however, for large enough Reynolds numbers, we typically include turbulence closure models for additional dissipation (Sec.~\ref{sec:turb}).

%The following table summarises the discretisations available in {\fontfamily{qcr}\selectfont  HORSES3D}:

\subsection{Implementation of alternate physics}
\label{sec:altenPhy}
The code was initially developed to solve the Navier-Stokes equations \eqref{eq:compressibleNScompact} in the form of a conservation law.
%
% \begin{equation}
%     \vec{u}_t  + \nabla\cdot\vec{\boldsymbol{F}}_e = \nabla\cdot\vec{\boldsymbol{F}}_v.
% \end{equation}
%
Most of the physical systems introduced in Sec.~\ref{sec:multiphysics} can be described within this framework with alternate definitions of the conserved variables $\boldsymbol{u}$, convective fluxes $\ssvec{F}_e$ and viscous fluxes $\ssvec{F}_v$. Of course, changing the definitions of these quantities affects the numerical part of the solver, e.g., the convective Riemann fluxes used. 
Among the exceptions, we have the multiphase system described in Sec.~\ref{sec:Mul}, which contains non-conservative terms that cannot be described within this framework, see details in Sec.~\ref{sec:Mul}.  

All the code that includes the definition of conserved variables, fluxes, and the numerics associated (e.g., numerical fluxes) is organised into a library called \emph{Physics}, so the set of equations that the code solves can be changed by modifying the Physics library. The alternate systems described in Sec.~\ref{sec:multiphysics} have been constructed in this way. At compilation time, different executables are generated (one for each different physics implemented in {\fontfamily{qcr}\selectfont  HORSES3D}). This framework allows for some flexibility, while keeping the efficiency of a one-physics-tailored code. 

\subsubsection{The problem file}

{\fontfamily{qcr}\selectfont  HORSES3D} incorporates a dynamic library that is compiled separately (after) from the main code compilation. This library includes several useful subroutines that enhance the interaction level of the user with the code and provides interfaces so that the core solver does not have to be altered for special cases. The main subroutines included are:
\begin{tcolorbox}[colback=black!5!white,colframe=gray!75!black]
\begin{itemize}
    \item{UserDefinedStartup:} a subroutine that is called before any other routines in the code. 
    \item{UserDefinedFinalSetup:} a subroutine that is called after the mesh is read to allow mesh-related initialisations or memory allocations.
    \item{UserDefinedInitialCondition:} a subroutine that is called to set the initial condition for the flow.
    \item{UserDefinedState/Grad/Neumann:} a set of subroutines that allows the definition of user-defined boundary conditions.
    \item{UserDefinedPeriodicOperation:} Called before/during/after/periodically at every time-step to allow periodic operations to be performed.
    \item{UserDefinedSourceTermNS:} Called to apply source terms to the system of equations.
    \item{UserDefinedFinalize:} Called after the solution is computed to allow, for example, error tests to be performed.
    \item{UserDefinedTermination:}  Called at the the end of the main driver after everything else is done.
\end{itemize}
\end{tcolorbox}

The problem file library allows the user to customise the behaviour of the solver for specific purposes. For instance, a pressure probe could be defined using specific interpolation operators
to extract the state at a specific location in the flow using 
UserDefinedFinalSetup. The necessary memory for this variable can be allocated in UserDefinedStartup, and post-processes (e.g., filtered) through UserDefinedPeriodicOperation. The final result could then be output in
the desired format using UserDefinedFinalize.
.

\subsection{Energy/Entropy stable versions}\label{sec:stable}
{\fontfamily{qcr}\selectfont  HORSES3D} can use Gauss-Lobatto (GL) quadrature points (in addition to Gauss) for the integral approximation, which, through the Summation-By-Parts Simultaneous-Approximation Term (SBP-SAT) property \cite{FISHER2013518,doi:10.1137/130932193}, allows one to recover the continuous property of energy/entropy conservation at a discrete level. This enhances the robustness of the simulations, especially in the presence of high numerical aliasing (e.g., under-resolved turbulence) \cite{DBLP:journals/jscic/ManzaneroRFVK18}. There is a plethora of work focused on the construction of entropy- and energy-stable schemes for the discontinuous Galerkin method \cite{FISHER2013518,doi:10.1137/130928650,10.1016/j.jcp.2016.09.013,10.1007/s10915-018-0702-1,MANZANERO2020109241,CHEN2017427,WINTERS20181,doi:10.1137/120890144} and the references therein. Useful reviews on energy/entropy stable schemes for discontinuous Galerkin schemes are \cite{10.1016/j.jcp.2016.09.013,CSIAM-AM-1-1}. 

\begin{tcolorbox}[colback=black!5!white,colframe=gray!75!black]
The following energy/entropy stable discretisations are available for the compressible Navier-Stokes solver in {\fontfamily{qcr}\selectfont  HORSES3D}:
Morinishi \cite{morinishi2010skew}, Ducros \cite{DUCROS2000114}, Kennedy-Gruber \cite{kennedy2008reduced}, Pirozzoli \cite{pirozzoli2010generalized}, Entropy conserving \cite{10.1016/j.jcp.2016.09.013}, Chandrasekar \cite{chandrashekar_2013,CHANDRASHEKAR2013527}.
\end{tcolorbox}

Recently, we have extended entropy-stable schemes to the Spalart-Allmaras RANS equations \cite{LODARES2022110998} and multiphase flows (i.e. Navier-Stokes coupled to Cahn-Hilliard with local $p-adaptation$ \cite{ntoukas2022entropy}), and these have been implemented in {\fontfamily{qcr}\selectfont  HORSES3D}.

\subsection{Spatial p-adaptation}
\label{sec:space_adapt}
The code allows the selection of a different polynomial order ($P$) for each element. 
Additionally, the polynomial order can be different in each spatial dimension ($P_x$, $P_y$, $P_z$), allowing anisotropic p-adaptation. 
Treatment of the geometrical representation of the p-non-conforming faces for general curvilinear hexahedral elements is performed following the rules presented in~\cite{kopriva2019free}.
The coupling between the faces of elements with different polynomial orders is performed using the mortar method \cite{kopriva2002computation}, which does not retain entropy stability. 
%TODO GERASIMOS: IS THAT RIGHT?  yes , iNS/CH solver is ES for non-conforming.IT WOULD BE NICE TO HAVE THIS FEATURE TOO FOR THE COMPRESSIBLE NS SYSTEM (needs modificaitons in the code).
Selection of the polynomial order inside each element can be performed before the simulation starts (along with the mesh generation process) or as an automated process for steady or unsteady simulations. The automatic adaptation process requires a sensor that specifies which elements should be refined/coarsened. Details of sensor computation and adaptation process are explained in the following sections.  

\subsubsection{Sensor for adaptation}

Our preferred adaptation algorithm included in {\fontfamily{qcr}\selectfont  HORSES3D} uses the truncation error, which is an efficient sensor for mesh adaptation \cite{roy2009strategies,fraysse2014quasi}. Similarly to adjoint mesh adaptation, truncation error adaptation targets the source of the discretisation error, rather than adapting polluted regions \cite{roy2009strategies}. Additionally, we have shown that the truncation error controls all functional errors (e.g., error in lift or drag prediction) \cite{kompenhans2016adaptation}.  

The truncation error can be estimated using the $\tau-$estimation technique. The estimation of the truncation error is cheap compared to adjoint techniques. The foundation of the $\tau-$ estimation for high-order was established in \cite{rubio2013estimation, rubio2015quasi}. This technique permits the estimation of the truncation error on a hierarchy of meshes and therefore allows the selection (in one step) of the correct polynomial order for each element in each spatial direction. $\tau-$estimation was recently improved \cite{rueda2019truncation} to be a byproduct of a multigrid cycle. 

%TODO GERASIMOS: remember to include the feature based (for example, simple gradient based) adaptation in the compressible NS. 
% Re: Include in the text or in {\fontfamily{qcr}\selectfont  HORSES3D}?

Automatic mesh adaptation can also be performed using a feature-based approach. This is applied to perform p-adaptation for the Cahn-Hilliard equation and for the incompressible Navier-Stokes/Cahn-Hilliard system (multiphase flows), see sections (Sec.~\ref{sec:iNS}) and (Sec.~\ref{sec:Mul}). The application of interest is to locally refine the region of the interface between the different phases, as it contains large gradients that need to be resolved to ensure the correctness of the solution.

\subsubsection{Adaptation process}

\begin{figure}[!htbp]
\includegraphics[width=0.7\textwidth]{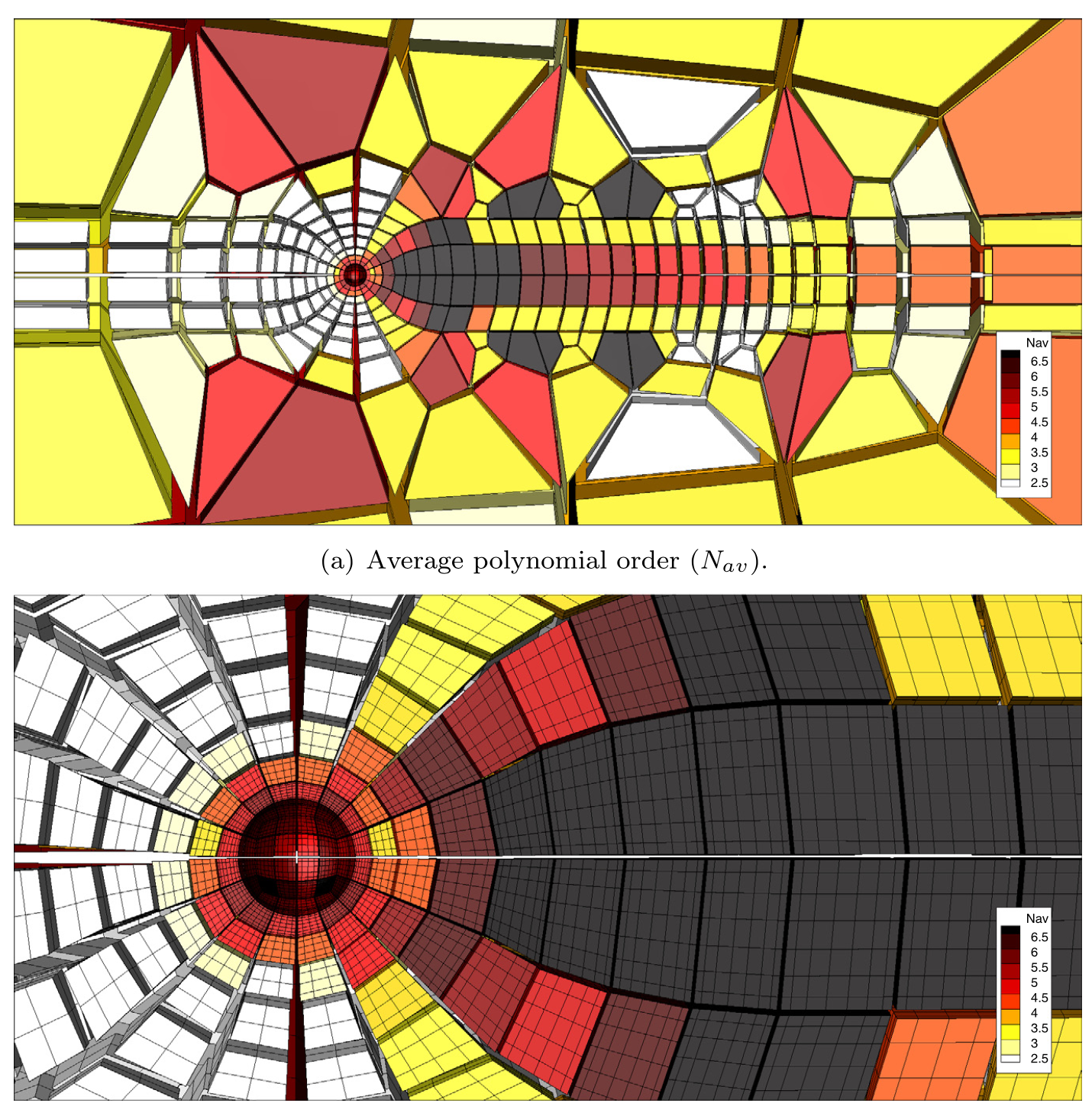}
\caption{Example of a p-adapted mesh for the flow around a sphere at Reynolds 200. The contours indicate the average polynomial order (Pav=(P1 +P2 +P3)/3). Taken from \cite{rueda2019p}.}
\label{fig:adaptation}
\end{figure}

%TRUNCATION ERROR adaptation
The truncation error adaptation algorithm performs mesh p-adaptation for steady and unsteady simulations. For steady flows, a simulation is conducted with a moderate polynomial order, $P$, and converged until the residual is smaller than the truncation error threshold set as the objective for the adapted mesh. Then, the truncation error is estimated with the $\tau$-estimation technique for all the polynomial orders coarser than the initial one. For every element, the lowest polynomial order (or polynomial order combination for anisotropic mesh adaptation) that meets the truncation error threshold is selected. If the polynomial orders from $1$ to $P$ do not meet the truncation error threshold, an extrapolation process is performed to select the correct polynomial order. In Fig.~\ref{fig:adaptation} an example of a p-adapted mesh based on truncation errors is shown. 
See \cite{kompenhans2016adaptation, kompenhans2016comparisons, rueda2019p} for more examples of truncation error-based adaptation for the compressible Navier-Stokes equations. For unsteady flows, the process is similar; however, no initial converged simulation is required, the truncation error is estimated and the mesh is adapted every prescribed number of iterations \cite{ntoukas2022entropy}. 

The adaptation process requires the user to select a target truncation error threshold for the adapted mesh. This threshold does not have any direct physical or engineering meaning and has been seen as a drawback of the methodology in the past. However, we have recently shown \cite{laskowski2022functional} how to link the truncation error threshold with any flow functional (e.g., lift, drag), bypassing this limitation .

\begin{figure}
\centering
\begin{subfigure}[t]{.45\textwidth}
  \centering
  \includegraphics[width = 1\textwidth]{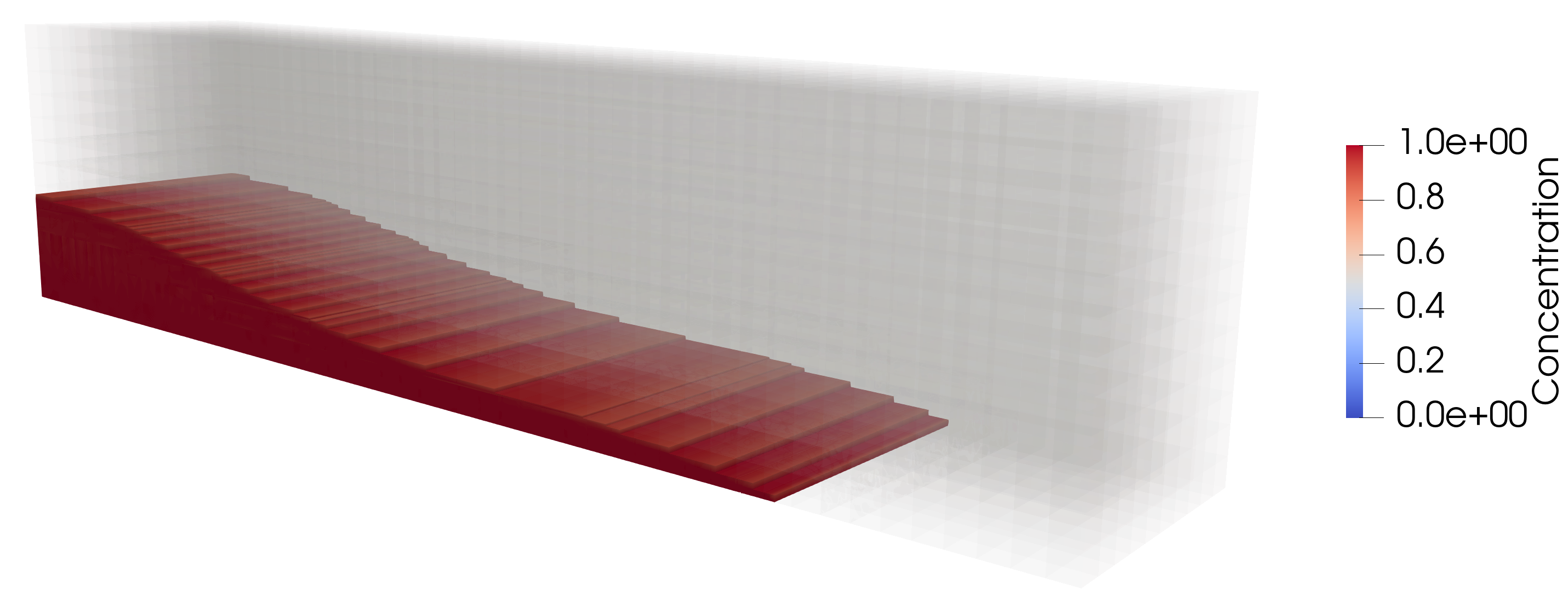}
  \caption{ $t=2$}
  \label{fig:BDInitialCondition}
\end{subfigure}%
\hfil
\begin{subfigure}[t]{.5\textwidth}
  \centering
  \includegraphics[width = 1\textwidth]{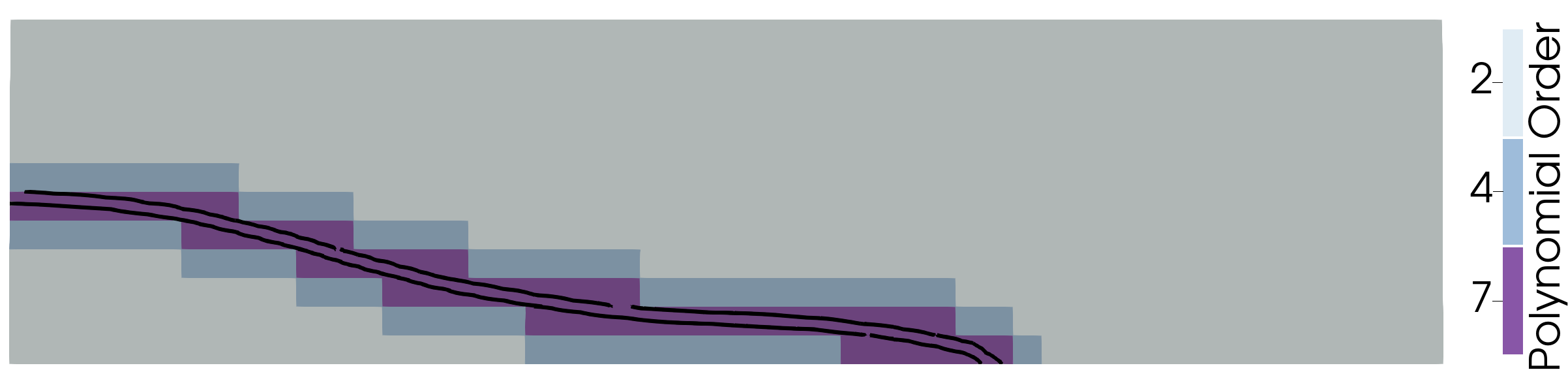}
  \caption{$t=2$}
  \label{fig:BD10}
\end{subfigure}
\caption{ Instances of the dam break test case showing the initial condition and the flow configuration at $t=2 \; s$. Taken from \cite{ntoukas2022entropy}.}
\label{fig:BD2Dvisualisation}
\end{figure}

% FEATURE BASED
The multiphase flow applications (Sec.~\ref{sec:Mul}) are of unsteady nature, and thus require the use of dynamic adaptation. In this case, it is crucial to have an adequate resolution in the region of the interface between different fluids. The tracking algorithm uses the value of the phase field parameter $\phi$ of the Cahn-Hilliard equation. The interface is defined as the region where $0.1\leq \phi \leq 0.9$. Thus, knowing the location of the interface throughout the simulation, we adapt the elements that contain part of the interface between different phases. This process, shown in Fig.~\ref{fig:BD2Dvisualisation}, is detailed in \cite{ntoukas2021free,ntoukas2022entropy}, as well as the method through which it can be automated.

\subsection{Time advancement}
\label{sec:time}

There are several temporal integration strategies in {\fontfamily{qcr}\selectfont  HORSES3D}, which are tailored for specific solvers. By default, the low-storage explicit Runge-Kutta of order 3 \cite{williamson1980low} is used as a time-stepping strategy for any given problem. This can be changed by the user to an explicit Runge-Kutta 4$^{th}$ order \cite{CarpenterKennedy1994} or to more sophisticated methods described below.

\begin{figure*}[htp]
    \centering
    \subfloat[Result for NACA0012 test case with 2$^{\text{nd}}$ order discretisation.]{
        \includegraphics[width=0.45\linewidth]{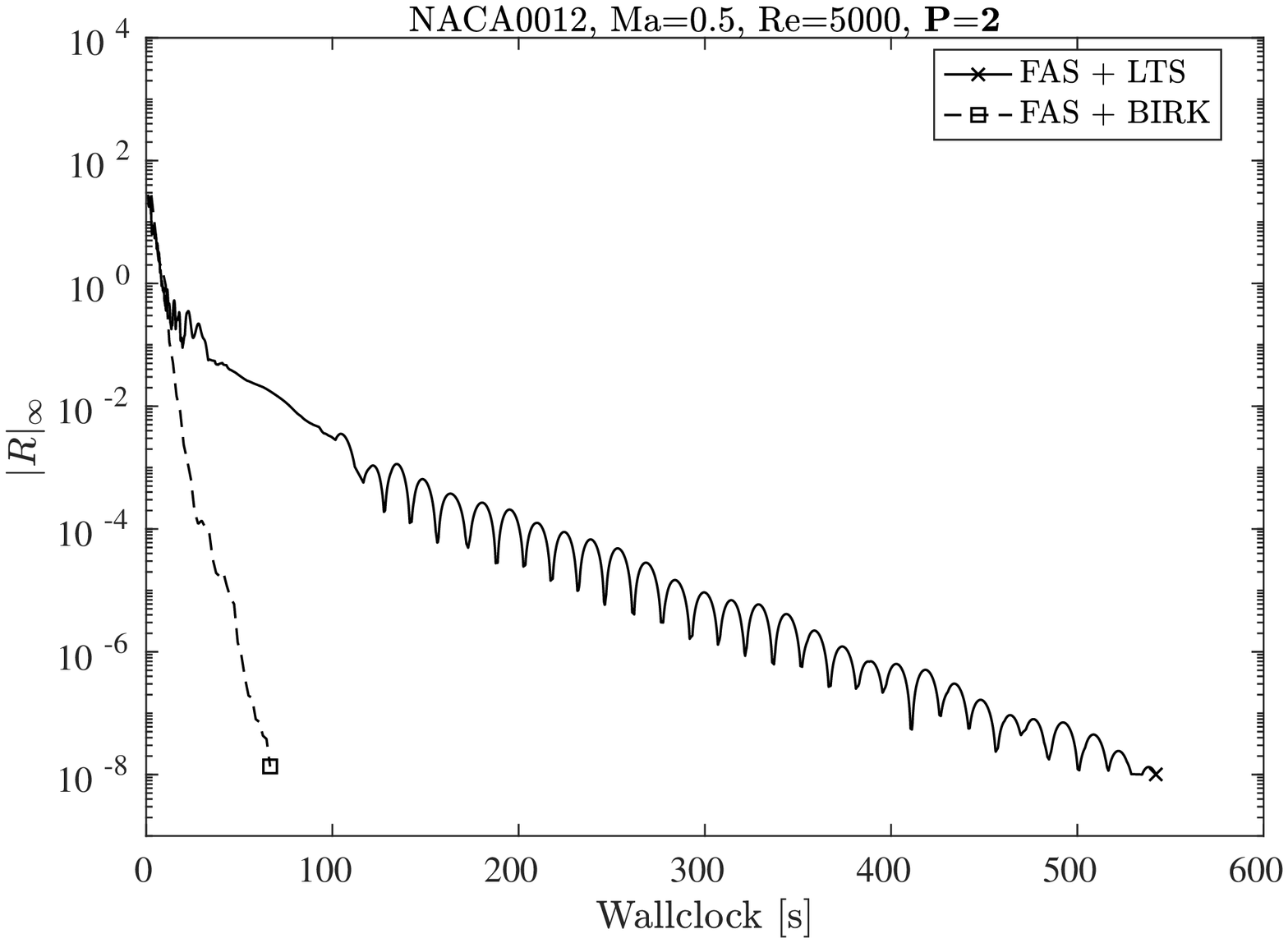}
        \label{fig:steady1}
    }
    \hfill
    \subfloat[Result for NACA0012 test case with 5$^{\text{th}}$ order discretisation.]{
        \includegraphics[width=0.45\linewidth]{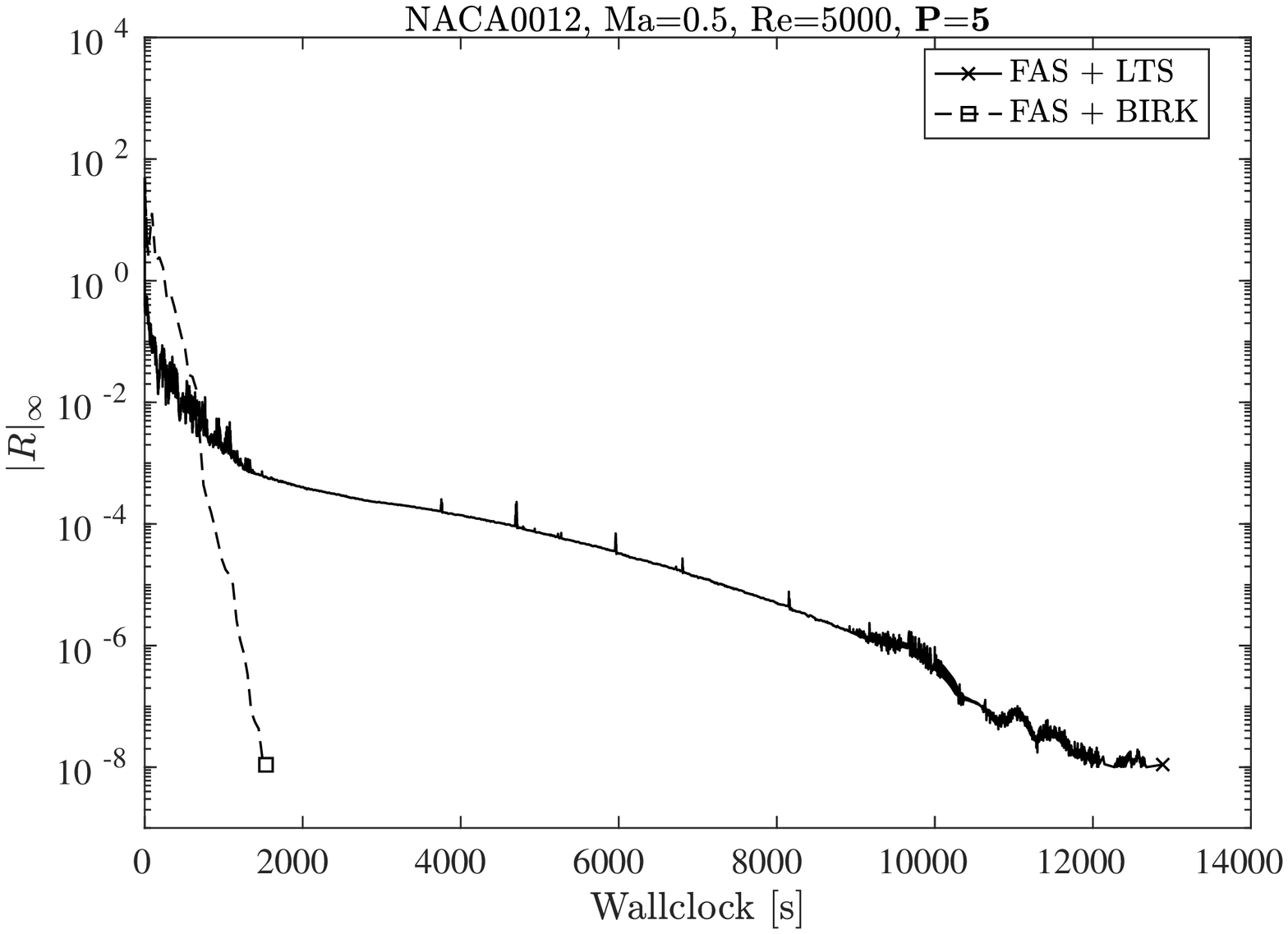}
        \label{fig:steady2}
    }
	\caption{Figures present convergence for a laminar test case (NACA0012 airfoil with $\text{Re}=5000$, $\text{Ma}=0.5$ and $\text{AoA}=0$) using FAS with explicit and semi-implicit convergence acceleration. The explicit FAS uses a five stage 2$^{\text{nd}}$ order Runge-Kutta smoother \cite{Vermeire2019} with local time-stepping. The semi-implicit FAS uses the BIRK smoother \cite{Ghidoni2014}. A low (left figure) and high (right figure) uniform polynomial discretisations are considered.}
    \label{fig:steadyresults}
\end{figure*}

\begin{figure}[!htbp]
\includegraphics[width=0.99\textwidth]{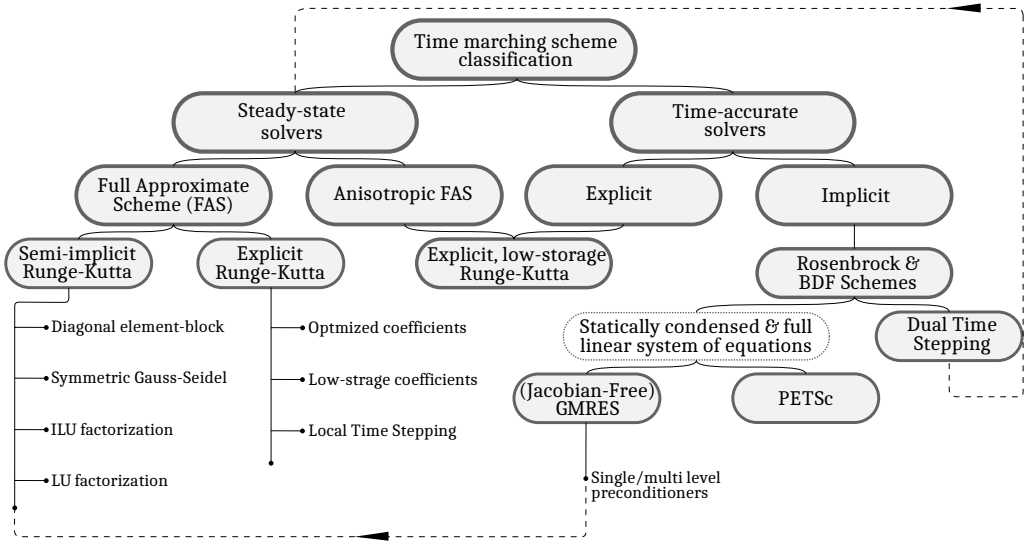}

\caption{Classification of time-marching schemes employed in {\fontfamily{qcr}\selectfont  HORSES3D}.}
\label{fig:solvers}
\end{figure}

For steady-state simulations, the non-linear p-multigrid method (also known as Full Approximate Scheme -- FAS) \cite{rueda2019p} is employed. FAS employs several explicit (local time-stepping, optimised Runge-Kutta coefficients for steady-state) \cite{Joshi2021,Vermeire2019} and semi-implicit (BIRK, LU-SGS, ILU(0)) \cite{Fidkowski2005,Parsani2010,Ghidoni2014} convergence acceleration strategies, presented in Fig.~\ref{fig:solvers}. An example of convergence for a laminar test case is presented in Fig.~\ref{fig:steadyresults}.

For time-accurate simulations (e.g., Large Eddy Simulations), implicit time integration can be used with Backward Difference Formulas (BDF) of order $1,...,5$, or 6$^{th}$ order Rosenbrock-type Runge-Kutta schemes \cite{Bassi2015}. The Jacobian can be computed analytically for the Navier-Stokes equations and also numerically when other equations need to be included (e.g., Spalart-Allmaras) using a coloring algorithm \cite{coleman1983estimation,gebremedhin2005color} to improve performance.

% \FloatBarrier

\subsection{Parallelisation and efficiency}
\label{sec:parallel_cost}
To reduce the computational cost, all techniques are parallelised for hybrid systems to run on multiprocessor CPUs combining MPI and OpenMP. Let us note that DG methods have compact stencils for unstructured meshes, and hence are highly parallelisable since each element has many structured internal degrees of freedom (related to the polynomial degree).

\begin{figure}[!htbp]
\includegraphics[width=\textwidth]{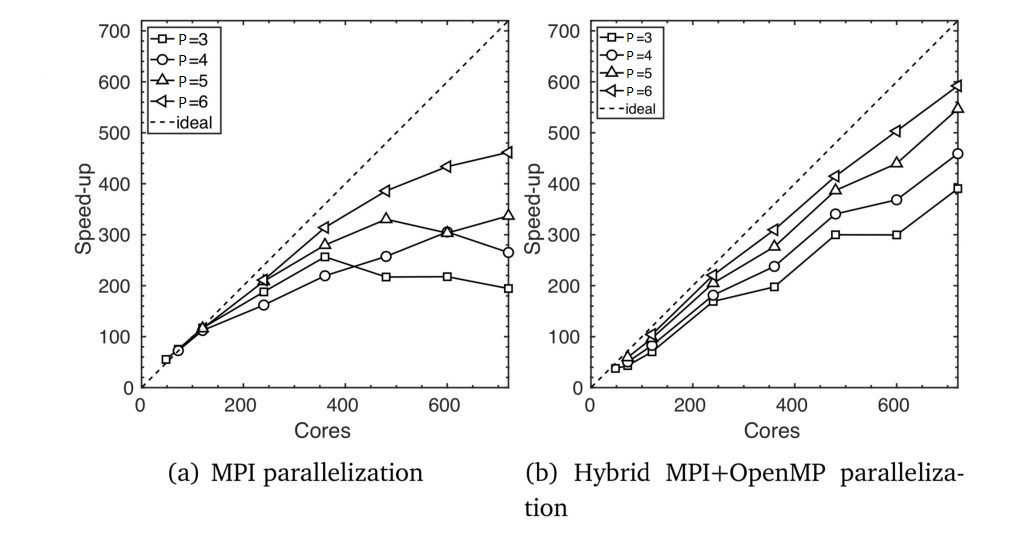}
\caption{{\fontfamily{qcr}\selectfont  HORSES3D} scalability solving the TGV problem with a $64^3$ mesh using a) MPI and b) MPI/OpenMP.}
\label{fig:parallel64}
\end{figure}

Numerical experiments for explicit time-stepping have been performed on the Taylor-Green Vortex (TGV) problem \cite{1937:Taylor}. The TGV problem has been widely used to report the subgrid-scale modelling capabilities of ILES approaches and discretisations \cite{WINTERS20181,2017:Moura,MANZANERO2020104440}. The results are shown in Fig. %\ref{fig:parallel32} and
\ref{fig:parallel64}.

We observe very good scalability when using high polynomial orders and MPI + OpenMP architectures up to a thousand cores; see Fig.~\ref{fig:parallel64}.b, where we get 86\% of the ideal performance for 600 cores and $P=6$.

%\begin{figure}[!htbp]
%\includegraphics[width=\textwidth]{Figures/parallel/parallel1.jpg}
%\caption{{\fontfamily{qcr}\selectfont  HORSES3D} strong scalability solving the TGV problem with a $32^3$ mesh using OpenMP.}
%\label{fig:parallel}
%\end{figure}

%\begin{figure}[!htbp]
%\includegraphics[width=\textwidth]{Figures/parallel/parallel2.jpg}
%\caption{{\fontfamily{qcr}\selectfont  HORSES3D} strong scalability solving the TGV problem with a $32^3$ mesh using MPI/OpenMP.}
%\label{fig:parallel32}
%\end{figure}

\section{Multiphysics}\label{sec:multiphysics}

In this section, we detail the various multiphysics applications included in {\fontfamily{qcr}\selectfont  HORSES3D}. Namely, we detail the turbulence closure models (RANS and LES), the incompressible flow formulation, multiphase modelling, shock capturing, flow-particle dynamics, aeroacoustics and immersed boundaries.

\subsection{Turbulence closure models}\label{sec:turb}
A particular challenge related to the numerical simulation of fluid flows relates to turbulence modelling. The main characteristic of turbulent flows is the presence of a broad spectrum, where a wide range of length and time scales coexist \cite{sagaut2006large}. It is the multiscale character (in time and space) of turbulent flows, which scales with the Reynolds number (i.e. ratio of inertia to viscous forces) that renders the numerical treatment of turbulence a difficult task. 

The Direct Numerical Simulation (DNS) approach solves the Navier-Stokes equations without additional modelling. Consequently, it requires a very high spatial and temporal resolution and is thus restricted to moderately low Reynolds. To compute large Reynolds number flows in complex configurations, it is necessary to introduce a certain degree of modelling to limit the cost. There are two main approaches: the Reynolds Averaged NS (RANS) approach and the Large Eddy Simulation (LES) technique. Both types are included in {\fontfamily{qcr}\selectfont  HORSES3D}.

\begin{figure*}[htp]
    \centering
            \subfloat[ Isometric view.]{
         \includegraphics[width=0.5\textwidth]{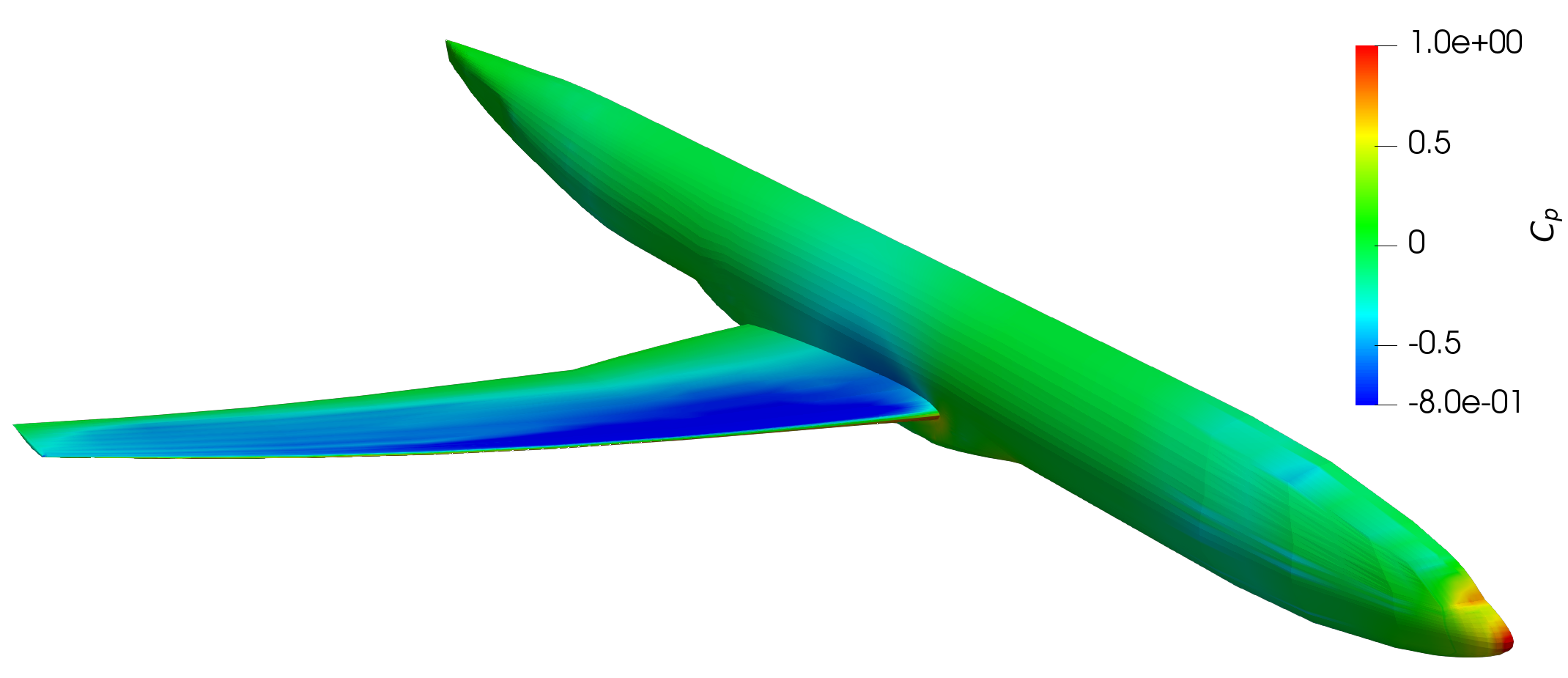}
    \label{fig:CRMisoview}
    }
    \subfloat[Top view. ]{
         \includegraphics[width=0.47\textwidth]{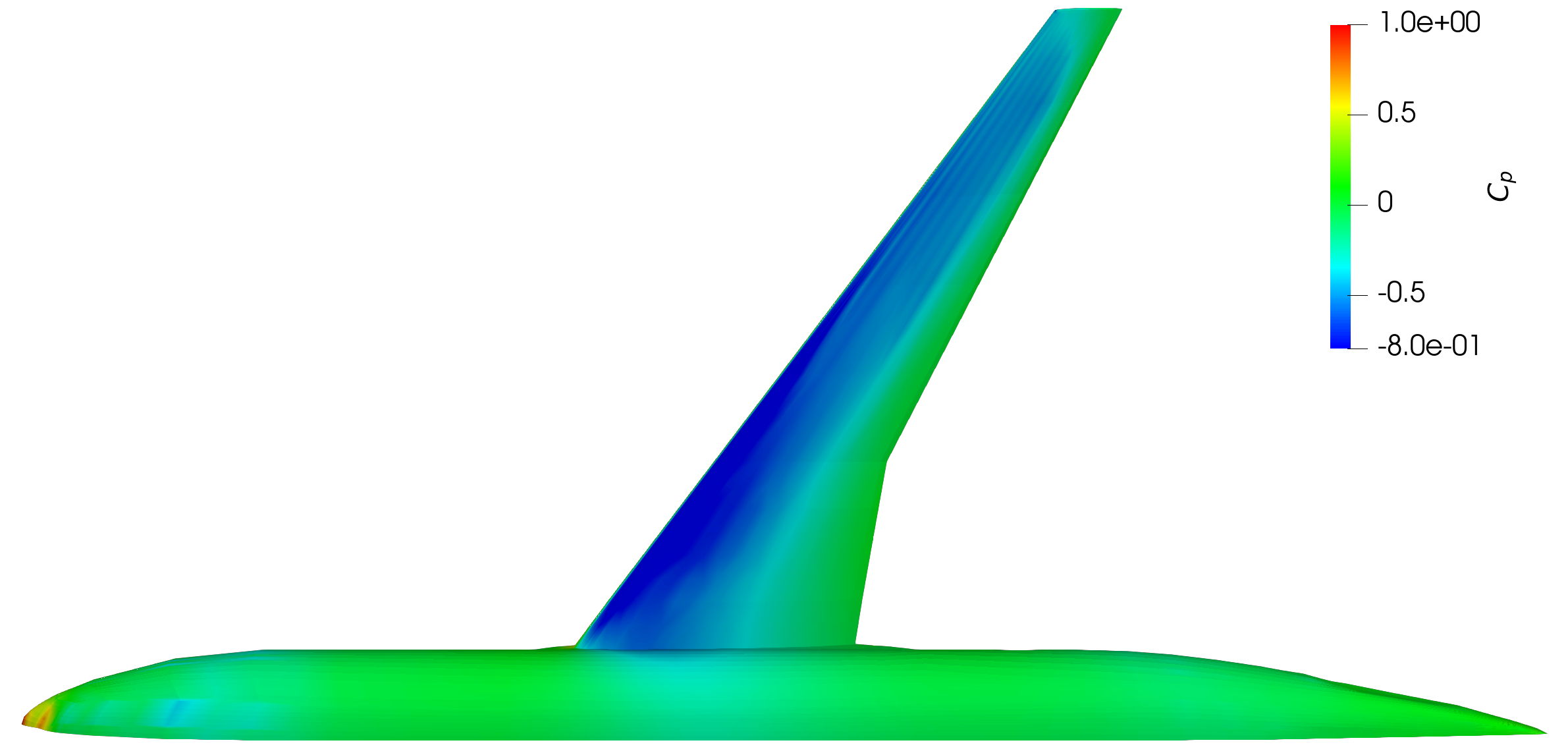}
        \label{fig:CRMtopview}
    }
	\caption{Pressure contour over the surface of the Common Research Model (CRM) for M=0.85, Re=$5\cdot10^6$ and AoA = $2.3^{\circ}$ using the {\fontfamily{qcr}\selectfont  HORSES3D} RANS module with split-form inviscid fluxes. The computations have been done on a $\mathrm{3^{rd}}$ order ($P=2$) hexahedral mesh from \cite{de2013robust}.} \label{fig:CRM}
\end{figure*}

The RANS solver uses the negative Spalart-Allmaras model from \cite{oliver2008high,Joshi2021}, while the LES models consider both explicit or implicit subgrid modelling. An example of the results derived using the {\fontfamily{qcr}\selectfont HORSES3D} framework is presented in Figure~\ref{fig:CRM}, where we consider the Common Research Model \cite{vassberg2008development}.  

The explicit LES approach uses spatial filtering where the large structures are resolved, reducing the modelling to the small unresolved turbulent structures (i.e. small eddies), which behave in an isotropic fashion. We have implementations for classic LES models, including Smagorinsky \cite{smagorinsky1963general,lilly1965computational}, Wale \cite{nicoud1999subgrid}, and Vreman \cite{vreman2004eddy}. In addition, we developed in \cite{MANZANERO2020104440} a SVV-Smagorinsky model where the Spectral Vanishing Viscosity technique was combined with a  Smagorinsky LES model to enhance the closure model. The model was capable of maintaining low dissipation levels for laminar flows, while providing the correct dissipation for all wave–number ranges in turbulent regimes. 

In recent years, Implicit Large Eddy Simulation (ILES) methods \cite{grinstein2007implicit} have gained popularity. This approach, increasingly popular when combined with high-order numerical methods (e.g., DG), uses the dissipation inherited from the numerical scheme to mimic turbulent effects of unresolved scales. DG methods have dissipation and dispersive errors, which are confined to high wavenumbers, and therefore introduce numerical dissipation primarily to under-resolved wave-lengths \cite{doi:10.1137/100807211}. This property makes them very suitable for ILES, since dissipation only acts at under-resolved scales, as explained in \cite{FerrerJCP,manzanero2018,Juan_vonNeumann,SOLANFUSTERO2021110246}.
Riemman solvers are the classic option to include numerical dissipation in DG schemes \cite{WorkshopDG,Gassner_ILES}, since they naturally arise when discretising the non-linear terms. A comparison of different fluxes for homogeneous turbulence can be found in \cite{FLAD2017782,MANZANERO2020104440}. A different option is to modify the viscous terms to enhance the approximation's dissipative properties. The modification has been proposed in \cite{FerrerJCP} using an increased penalty parameter (compared to the minimum required to ensure coercivity of the scheme) when discretising the viscous terms using an interior penalty formulation.

\begin{figure}[!htbp]
\includegraphics[width=0.99\textwidth]{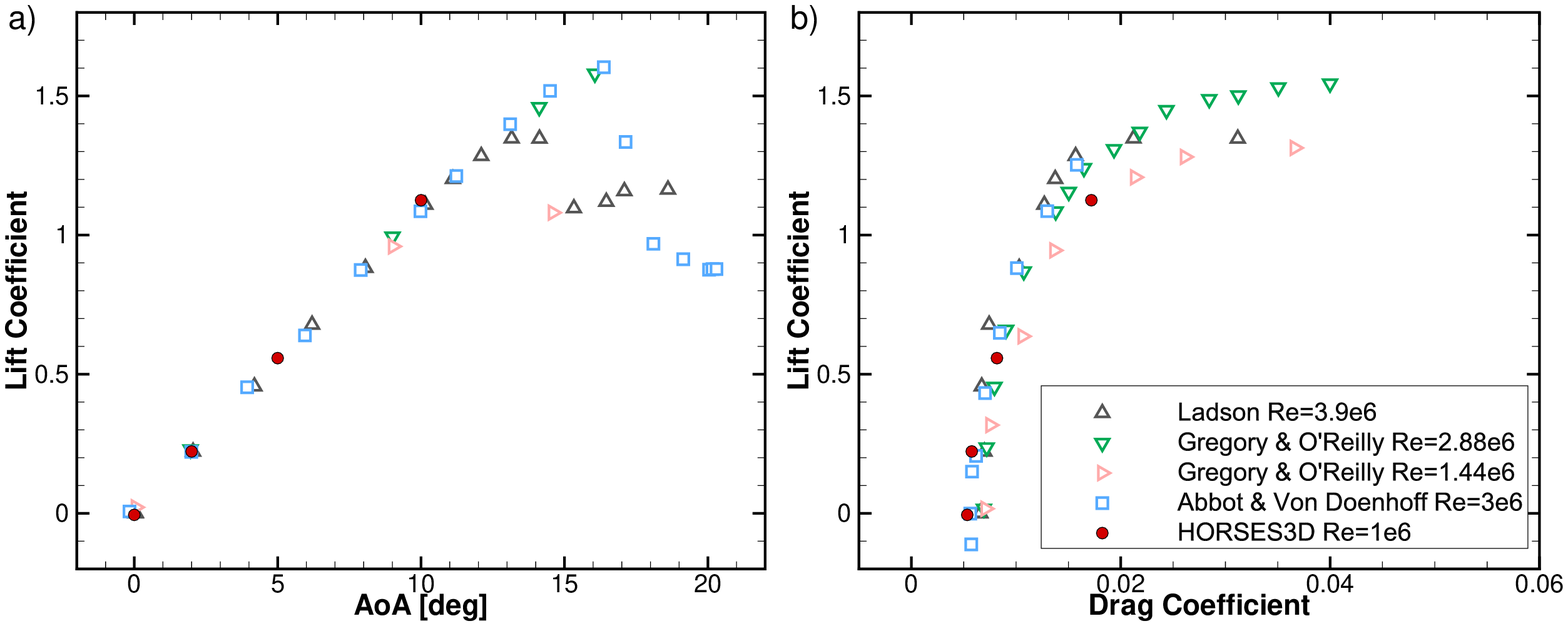}

\caption{Lift and drag for a NACA 0012 at Re=1e6;  comparisons between {\fontfamily{qcr}\selectfont  HORSES3D} (ILES version) with experimental data. Adapted from \cite{10.1007/978-3-030-39647-3_38}}
\label{fig:lift_drag}
\end{figure}

When combined with high-order DG methods, ILES provide valuable results even for challenging problems such as the prediction of turbulent quantities, detached flows and laminar-turbulent transition on airfoils, see, for example, \cite{Uranga,FerrerJCP,10.1007/978-3-030-39647-3_38,grinstein2007implicit,MANZANERO2020104440,FERNANDEZ2017308,FERRER2019104239}. Finally, a promising approach is to combine robust energy-stable methods (e.g., skew symmetric variants) with ILES techniques (e.g., Spectral Vanishing Viscosity) to provide accurate simulation at high Reynolds numbers, see \cite{manzanero2018,MANZANERO2020104440}. 
As an example of ILES performance, we calculated the averaged values for the lift and drag on a NACA0012, using a hexahedral mesh with 18000 elements, with $P=4$ (2.2 million degrees of freedom).
Fig. \ref{fig:lift_drag} compares the aerodynamic coefficients with the experimental data for various angles of attack and the two solvers. Fig. \ref{fig:lift_drag} [Left] shows the lift coefficient against $\text{AoA}$ and Fig. \ref{fig:lift_drag} [Right] depicts the lift-drag polar for $\text{Re}=1\times10^6$. We observe very good agreement with the experimental data.
% \FloatBarrier

\begin{tcolorbox}[colback=black!5!white,colframe=gray!75!black]
The following turbulence models are available in {\fontfamily{qcr}\selectfont  HORSES3D}:
RANS: Spalart-Allmaras; explicit LES: Smagorinsky, Wale, Vreman, SVV-Smagorinsky; implicit LES: energy/entropy stable formulations and dissipative fluxes.
\end{tcolorbox}

\subsection{Incompressible Flows}
\label{sec:iNS}
Among the various incompressible Navier-Stokes models, {\fontfamily{qcr}\selectfont  HORSES3D} uses an artificial compressibility formulation \cite{shen1997pseudo}, which converts the elliptic problem into a hyperbolic system of equations, at the expense of a non divergence–free velocity field. However, it allows one to avoid constructing an approximation that satisfies the inf–sup condition, see \cite{Sherwin,ferrer_moxey_willden_sherwin_2014} and the references therein. The artificial compressibility model is commonly combined with dual time–stepping, which solves an inner pseudo-time-step loop until velocity divergence errors are lower than a selected threshold and then performs the physical time marching. This methodology is well suited for use as a fluid flow engine for interface–tracking multiphase flow models, as it allows the density to vary spatially $\rho\left( \vec{x},t \right)$. The artificial compressibility system of equations is
\begin{equation}
  \rho_{t} + \vec{\nabla}\cdot\left(\rho\vec{u}\right)=0 ,
  \label{eq:incomressible:continuity}
\end{equation}
\begin{equation}
\left(\rho\vec{u}\right)_t+\vec{\nabla}\cdot\left(\rho \vec{u}\vec{u}\right) = -\vec{\nabla}p + \vec{\nabla}\cdot\left(\frac{1}{\mathrm{Re}}\left(\vec{\nabla}\vec{u} + 
\vec{\nabla}\vec{u}^{T}\right)\right)+\frac{1}{\mathrm{Fr^{2}}}\rho\vec{e}_{g},
\label{eq:incomressible:momentum}
\end{equation}
\begin{equation}
  p_t + c_0^2 \vec{\nabla}\cdot\vec{u} = 0,
  \label{eq:incomressible:ACM}
\end{equation}

where Re is the Reynolds number, Fr is the Froude number, $c_0^2$ is the artificial speed of sound, and $\vec{e}_g$ is the gravitational acceleration.
%To increase the robustness of the solver, particularly for under-resolved flows, an entropy-stable discretisation has been implemented. This can be achieved through a split-form nodal DG that uses the summation-by-parts/simultaneous-approximation-term (SBP-SAT) property. Different averaging options are considered for the two-point inviscid fluxes in \cite{Juan}. The detailed stability analysis is presented in \cite{Juan} along with test cases for variable density flows such as the Rayleigh-Taylor instability. 

\subsection{Multiphase}
\label{sec:Mul}
Phase field models describe the phase separation dynamics of two immiscible liquids by minimising a chosen free–energy. For an arbitrary free–energy function, it is possible to construct different phase field models. Among the most popular, one can find the Cahn–Hilliard \cite{cahn1958free} and the Allen–Cahn \cite{allen1972ground} models. The popularity of the first, despite being a fourth-order operator in space, comes from its ability to conserve the phases. {\fontfamily{qcr}\selectfont  HORSES3D} permits the solution of the Cahn-Hilliard equation. For a detailed explanation of the model, see \cite{manzanero2020free}.

The multiphase flow solver implemented in {\fontfamily{qcr}\selectfont  HORSES3D} is constructed by a combination of the diffuse interface model of Cahn–Hilliard \cite{cahn1958free} with the incompressible Navier–Stokes equations with variable density and artificial compressibility \cite{shen1997pseudo}. This model is entropy stable and guarantees phase conservation with an accurate representation of surface tension effects. For a detailed explanation of the model, see \cite{MANZANERO2020109241}. The modified entropy-stable version approximates 
\begin{equation}
c_t + \vec{\nabla}\cdot\left(c\vec{u}\right) = M_0 \vec{\nabla}^2 \mu,
\label{eq:governing:cahn--hilliard}
\end{equation}
\begin{equation}
\sqrt{\rho}\left(\sqrt{\rho}\vec{u}\right)_t+\vec{\nabla}\cdot\left(\frac{1}{2}\rho \vec{u}\vec{u}\right) 
+\frac{1}{2}\rho\vec{u}\cdot\vec{\nabla}\vec{u}+c\vec{\nabla}\mu
= -\vec{\nabla}p + \vec{\nabla}\cdot\left(\eta\left(\vec{\nabla}\vec{u} + 
\vec{\nabla}\vec{u}^{T}\right)\right)+\rho\vec{g},
\label{eq:governing:momentum-skewsymmetric-sqrtRho}
\end{equation}
\begin{equation}
  p_t + \rho_0 c_0^2 \vec{\nabla}\cdot\vec{u} = 0,
  \label{eq:governing:ACM}
\end{equation}
where $c$ is the phase field parameter, $M_0$ is the mobility, $\mu$ is the chemical potential, $\eta$ is the viscosity and $c_0$ is the artificial speed of sound. 
In \cite{MANZANERO2020109241}, we provide a comparison between standard and entropy-stable discretisations. Through several numerical experiments, the superior robustness characteristics of the entropy-stable scheme are highlighted in comparison to the standard scheme. Among the test cases considered, there is a two--phase flow of oil and water within a pipe \cite{MANZANERO2020109241} and a 3D dam break test case \cite{ntoukas2022entropy}. 

\begin{figure}
    \centering
    \includegraphics[width=0.9\linewidth]{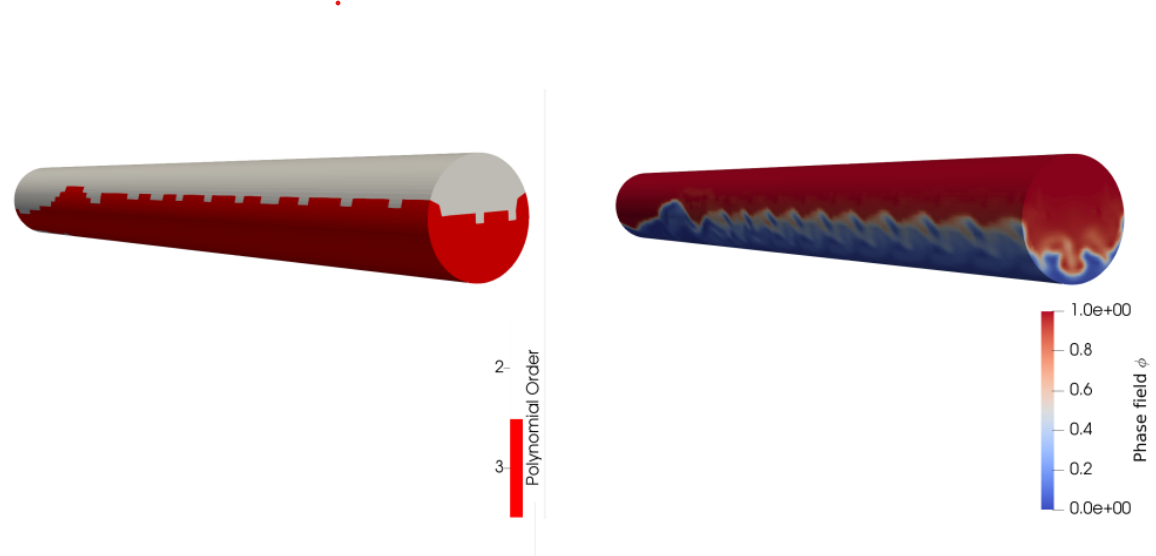}
    \caption{Annular flow test case of water and oil. On the left, the polynomial order distribution is presented for the p-adaptive scheme based on the interface tracking refinement. On the right, the flow structure is presented for the developed flow. Taken from \cite{ntoukas2022entropy}. }
    \label{fig:multi:pipe}
\end{figure}

For phase field modelling, it is important to adequately resolve the region of the interface between different phases, as it contains large gradients. In {\fontfamily{qcr}\selectfont  HORSES3D}, the local refinement around the interface is carried out through p-adaptation. The original scheme presented in \cite{MANZANERO2020109241} has been extended to support p-non-conforming elements and has been modified accordingly so that it is entropy-stable even for p-non-conforming elements, see Fig.~\ref{fig:multi:pipe}. The detailed analysis and modifications to the original scheme are presented in \cite{ntoukas2022entropy}. 

\subsection{Shock capturing}
\label{sec:Shock}

Supersonic flows present additional difficulties to high-order spectral elements, where the favourable accuracy and convergence properties of these schemes can be damaged by the Gibbs phenomenon near shocks. {\fontfamily{qcr}\selectfont  HORSES3D} implements a special type of artificial viscosity to handle shocks. In the case of the Navier-Stokes equations, discontinuities due to compressibility effects have an actual width proportional to the viscosity,~$\mu$ \cite{vonMises1950}. This motivates the introduction of an artificial viscous term ~$\ssvec{F}_a$ into the equations to control the effective~$\mu$ at each element to smooth the solution.

{\fontfamily{qcr}\selectfont  HORSES3D} includes two options for $\ssvec{F}_a$: the viscous flux of the Navier-Stokes equations and the Guermond-Popov entropy-stable flux \cite{Guermond2014}, which also includes a regularisation term in the density equation. Regardless of the selected viscous flux, and to minimise artificial viscosity in smooth regions, we introduce an Spectral Vanishing Viscosity (SVV) formulation and apply a frequency-dependent filter kernel~$\mathcal{F}$ to these fluxes, so that the additional term has the form~$\ssvec{F}_a^{\mathcal{F}} = \mathcal{F}\star\ssvec{F}_a$. We implement the filter kernels of \cite{tadmor_89,Maday1993,Moura2}. Since the code is a nodal DGSEM framework, we transform the fluxes to the Legendre modal space,~$\{L_i\}$, where the application of this filter is simply an element-wise product of the fluxes and the kernel, and then we return to the nodal space to recover the filtered flux.

The stability of this new term is analysed in terms of the evolution of the entropy. Physically plausible solutions of the Navier-Stokes equations must agree with elemental laws such as the conservation of the energy or a non-decreasing entropy. We impose this by defining a mathematical entropy function that represents the underlying law and requiring it to be bounded. In this case, the mathematical entropy,~$S$, is defined in terms of the non-decreasing physical entropy as
\begin{equation}
    S = -\rho s, \quad s = \ln p - \gamma \ln \rho.
\label{eq:shock:entropy}
\end{equation}
We also introduce the concept of entropy variables~$\boldsymbol{w}$, which are computed from~\eqref{eq:shock:entropy} as
\begin{equation}
    \boldsymbol{w} = \frac{\partial S}{\partial \boldsymbol{u}}, \quad
    \boldsymbol{u} = (\rho, \rho\vec{v}, \rho e)^T.
\end{equation}

To ensure the stability of the new term, a detailed semi-discrete entropy analysis (see \cite{Friedrichs1971,Tadmor1986,LODARES2022110998}) shows that filtering has to be applied in a very specific form. Expressing the baseline flux as
\begin{equation}
    \ssvec{F}_a = \mathcal{B}\ssvec{G}, \quad
    \mathcal{B} = \mathcal{L}^T\mathcal{D}\mathcal{L}, \quad
    \ssvec{G} = \nabla\boldsymbol{w},
\label{eq:shock:fluxmatrix}
\end{equation}
where~$\mathcal{B}$ is a symmetric matrix, and~$\mathcal{L}$ and~$\mathcal{D}$ are the matrices of its $\textbf{L}^T\textbf{D}\textbf{L}$ decomposition, the filtered flux,~$\ssvec{F}_a^{\mathcal{F}}$, is obtained by filtering only one half of~\eqref{eq:shock:fluxmatrix},
\begin{equation}
    \ssvec{F}_a^{\mathcal{F}} = \frac{1}{\sqrt{J}}\mathcal{L}^T\sqrt{\mathcal{D}}\mathcal{F}\star
        \left(\sqrt{J\mathcal{D}}\mathcal{L}\ssvec{G}\right),
\label{eq:shock:filteredflux}
\end{equation}
where $J$ is the Jacobian of the transformation from the reference element to the physical elements of the mesh, see Sec.~\ref{sec:dg}. A detailed explanation of this derivation can be found in \cite{Mateo2022}.

\begin{figure}
    \centering
    \includegraphics[width=0.8\linewidth]{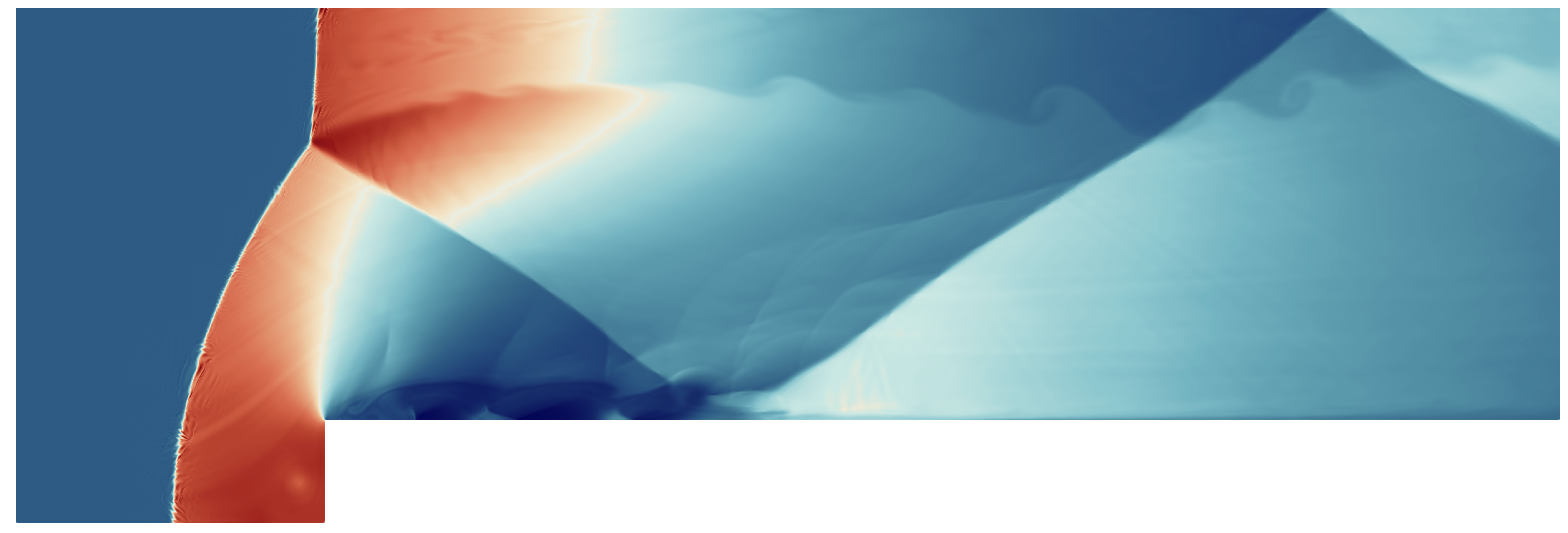}
    \caption{Facing forward step case at~$\tau=10$ with an incoming flow at Ma=3 and 3653 elements with polynomial degree~$P=7$. Adapted from \cite{Mateo2022}.}
    \label{fig:shock:ffs}
\end{figure}

The SVV approach has proven to be useful for the simulation of turbulent and supersonic flows (see Figure~\ref{fig:shock:ffs}); however, the treatment of discontinuities requires a higher amount of dissipation that can be achieved only by adding non-filtered dissipation in troubled regions. Since this is only required near discontinuities, we use a sensor,~$s_{\rho}$, based on the average value of the density gradient to detect troubled elements and modify the filter there,
\begin{equation}
    s_{\rho} = \int_{el} J\lvert\nabla_{\xi}\rho\rvert^2.
\end{equation}

%\begin{figure}
%    \centering
%    \includegraphics[width=0.5\linewidth]{Figures/shock/TGV.png}
%    \caption{Q-contours of the Taylor-Green vortex at~$\tau=5$.}
%    \label{fig:shock:tgv}
%\end{figure}

\subsection{Particles}\label{sec:particles}

{\fontfamily{qcr}\selectfont  HORSES3D} includes a two-way coupled Lagrangian solver. The model presented here is based on \cite{frankel2016settling}. 

Particles are tracked along their trajectories, according to the simplified particle equation of motion, where only contributions from Stokes drag and gravity are retained,
\begin{equation}
\label{eq:part_motion}
\frac{d y_i}{dt} = u_i, \quad \frac{d u_i}{dt} = \frac{v_i - u_i}{\tau_p} + g_i,
\end{equation}
where $u_i$ and $y_i$ are the \emph{ith} components of velocity and position of the particle, respectively. Furthermore, $v_i$ accounts for the continuous velocity of the fluid at the position of the particle.  We consider spherical Stokesian particles, so their mass and aerodynamic response time are $m_p = \rho_p \pi D_p^3/6$ and $\tau_p = \rho_p D_p^2 / 18\mu $, respectively, $\rho_p$ being the particle density and $D_p$ the particle diameter. 

Each particle is considered to be subject to a radiative heat flux $I_o$. The carrier phase is transparent to radiation, whereas the incident radiative flux on each particle is completely absorbed. Because we focus on relatively small volume fractions, the fluid-particle medium is considered to be optically thin. Under these hypotheses, the direction of the radiation is inconsequential, and each particle receives the same radiative heat flux, and its temperature $T_p$ is governed by
\begin{equation}
\label{eq:part_energy}
\frac{d}{dt} (m_p c_{V,p} T_p) = \frac{\pi D_p^2}{4} I_o - \pi D_p^2 h (T_p-T),
\end{equation}
where $c_{V,p}$ is the specific heat of the particle, which is assumed to be constant with respect to temperature. $T_p$ is the particle temperature and $h$ is the convective heat transfer coefficient, which for a Stokesian particle can be calculated from the Nusselt number $Nu = hD_p/k = 2$.

In practical simulations, integrating the trajectory of every particle is too expensive. Therefore, particles are agglomerated into parcels, each of them accounting for many particles with the same physical properties, position, velocity, and temperature. The evolution of the parcels is tracked with the same set of equations presented for the particles.

The two-way coupling means that fluid flow is modified because of the presence of particles. Therefore, the Navier-Stokes equations \eqref{eq:compressibleNScompact} are enriched with the following source terms:
\begin{equation}
\boldsymbol{S} = \beta\left[\begin{array}{c} 0 \\
                                                                       \sum_{n=1}^{N_p} \frac{m_p}{\tau_p} (u_{1,n}-v_1)\delta(\mathbf{x} - \mathbf{y}_n)  \\
                                                                       \sum_{n=1}^{N_p} \frac{m_p}{\tau_p} (u_{2,n}-v_2)\delta(\mathbf{x} - \mathbf{y}_n) \\
                                                                      \sum_{n=1}^{N_p} \frac{m_p}{\tau_p} (u_{3,n}-v_3)\delta(\mathbf{x} - \mathbf{y}_n) \\
                                                                      \sum_{n=1}^{N_p} \pi D_p^2 h (T_{p,n} - T) \delta(\mathbf{x} - \mathbf{y}_n )
\end{array}\right],
\label{eq:particlessource}
\end{equation}
where $\delta$ is the Dirac delta function, $N_p$ is the number of parcels, $\beta$ is the number of particles per parcel and $u_{i,n}$, $\mathbf{y}_{i,n}$, $T_{p,n}$ are the velocity, spatial coordinates, and temperature of the parcel \emph{nth}.

The ordinary differential equation (ODE) system \eqref{eq:part_motion}, \eqref{eq:part_energy} that controls the evolution of parcels is integrated in time using a low storage Runge-Kutta 3 method \cite{williamson1980low}. Integration is separated from the integration of the Navier-Stokes equations.
The properties of the fluid (velocity, viscosity, and temperature) at the position of the parcel are evaluated using the interpolation routines included in {\fontfamily{qcr}\selectfont  HORSES3D}. These routines work at the element level, so a routine to find the element to which the parcel belongs in general unstructured meshes is also included. 
The Dirac delta, $\delta(\mathbf{x} - \mathbf{y}_n)$, which appears in the source term, see \eqref{eq:particlessource}, can be dealt with in a simple way in the DG setting, since it can be integrated exactly in the weak form. 

\begin{figure}[!htbp]
\includegraphics[width=0.7\textwidth]{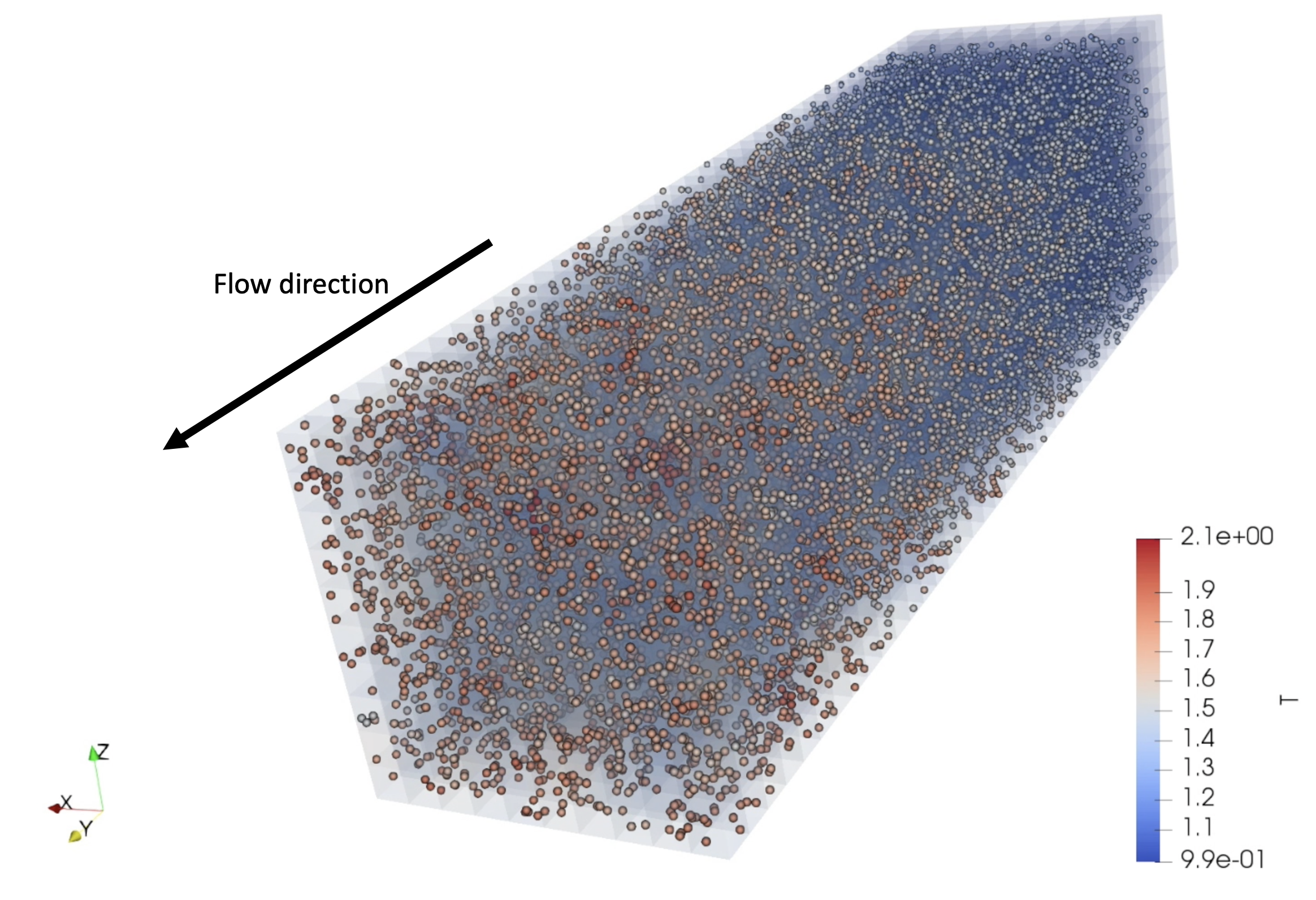}
\caption{Visualisation of a particle-laden flow in a channel subject to radiation. Particles and gas are injected into the channel at the same temperature. The temperatures of the particles increase in the direction of the flow as a result of the effect of radiation. Details of the test case are given in \cite{pouransari2017effects}. }
\label{fig:particles3d}
\end{figure}

\begin{figure}[!htbp]
\includegraphics[width=0.7\textwidth]{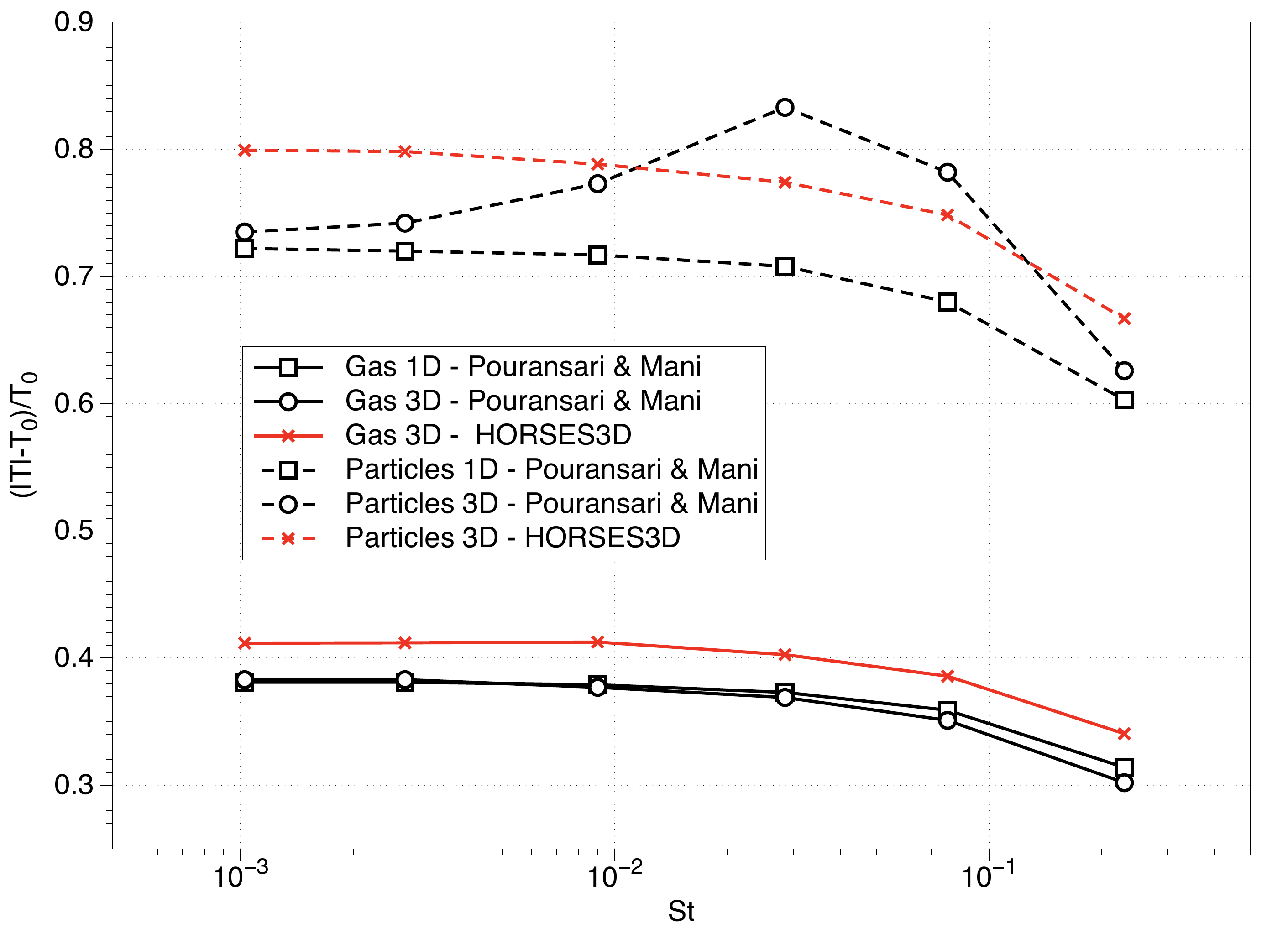}
\caption{Mean temperatures of gas and particles at the outlet of a channel subject to radiation. $T_0=300$ K is the reference temperature, while $St$ is the Stokes number. Details of the test case are given in \cite{pouransari2017effects}. }
\label{fig:particles}
\end{figure}

The coupled-particle model has been used to analyse a particle-laden flow in a channel subject to radiation, which is representative of particle solar radiation systems. Details of the test case can be found in \cite{pouransari2017effects}. The results obtained with {\fontfamily{qcr}\selectfont  HORSES3D} were obtained with a Cartesian mesh of 4000 elements and a number of particles per parcel $\beta=160$, see Fig. \ref{fig:particles3d}. The results of {\fontfamily{qcr}\selectfont  HORSES3D} are compared in Fig. \ref{fig:particles} with the one-dimensional and three-dimensional simulations performed in \cite{pouransari2017effects}. 
The results of {\fontfamily{qcr}\selectfont  HORSES3D} do not exactly match the three-dimensional simulations of \cite{pouransari2017effects} as the effect of preferential concentration at the inlet was not considered. Besides {\fontfamily{qcr}\selectfont  HORSES3D} uses parcels (agglomerations of particles). This artificial concentration of particles creates local peaks in temperature, which results in an increase of the average temperature.

\subsection{Computational Aeroacoustics}\label{sec:acoustic}

Computational aeroacoustics (CAA) can be divided into two categories. The first one, noise generation, is characterised by turbulence and non-linear interactions near solid surfaces. The second category is noise propagation and can be simplified to a linear propagation of the sound waves. %Each region has a different characteristic length of scale (and time scale)~\cite{Beck2018}.

An approach in CAA is to solve the near and far fields together, using the compressible Navier–Stokes equations on the whole domain, obtaining all the physical effects on the solution without needing external models. This approach is called Direct Noise Computation (DNC), and can deal with complex acoustics, such as acoustic feedback. The compressible Navier–Stokes of {\fontfamily{qcr}\selectfont  HORSES3D} is spectrally accurate and therefore well suited to provide accurate acoustics (e.g., pressure fluctuations).

The other possibility for CAA is to use acoustic analogies, which decouple the near and far fields into two computational regions, the so-called hybrid approach. Here, the CAA solution becomes a two-step process, where the near field is computed first and independently of the propagation~\cite{Lele2014} using the {\fontfamily{qcr}\selectfont  HORSES3D} compressible Navier–Stokes solver. Subsequently, equivalent acoustic sources can be extracted from the well-resolved near-field region, to feed the propagation problem, which is solved in the far field. This decoupling enables faster computations since the near-field that requires highly accurate flow computations can be performed separately from the acoustic propagation step. 

{\fontfamily{qcr}\selectfont  HORSES3D} can compute aeroacustics using the direct approach and also acoustic analogies. In particular, the 
 Ffowcs-Williams and Hawkings (FWH) aeroacoustic analogy~\cite{Ffowcs1969} is implemented. FWH rearranges the compressible NS equations into an inhomogeneous wave equation. The solution of this wave equation can be obtained in an integral representation (computed numerically).
The FWH uses a fictitious surface located in the near-field, where the Navier–Stokes solution is gathered and used to calculate the integral solution. This surface can coincide with the solid boundary of the body (having a solid surface) or be located at an arbitrary position in the flow field (being a permeable surface). {\fontfamily{qcr}\selectfont  HORSES3D} allows both possibilities.

%The implementation of FWH done in {\fontfamily{qcr}\selectfont  HORSES3D} has been done for cases similar to a wind tunnel scenario, which is when a stationary (permeable of solid) surface ($v_i=0$) is subjected to a constant flow in the positive streamwise direction ($U_{0i} = Constant \neq 0$), for any angle of attack.
Our implementation follows the formulation by Najafi~\cite{Najafi_Yazdi_2010}, and Ghorbanias~\cite{Ghorbaniasl_2012}, where the equations are recast as a convective wave equation. This formulation clearly separates the flow velocity fluctuations ($v_i$), surface velocities ($v_i$), and the free-stream velocity ($U_{0i}$) as follows:

\begin{equation} \label{conv_wave_eqn}
\begin{split}
    & \bigg( \frac{\partial^2}{\partial t^2} + U_{0i} U_{0j} \frac{\partial^2}{\partial x_i\partial x_j} + 2U_{0i} \frac{\partial^2}{\partial x_i\partial t} - c_0^2\frac{\partial^2}{\partial x_i\partial x_i} \bigg)(H(f)\rho ') = \\
     & \bigg( \frac{\partial}{\partial t} + U_{0i}\frac{\partial}{\partial x_i} \bigg) \big(Q_i n_i \delta (f) \big) - \frac{\partial}{\partial x_i} \big( L_{ij} n_j \delta (f) \big) + \frac{\partial^2}{\partial x_i\partial x_j}\big( T_{ij} H(f) \big),
\end{split}
\end{equation}
where the zero subscript refers to the undisturbed or free-stream flow, $\rho '$ is the density fluctuation, $Q_i$, $L_{ij}$, $T_{ij}$ are the loading, thickness, and quadrupole terms, respectively, $n_i$ is the $i^{th}$ component of the normal at the surface, $H(f)$ is the Heaviside function and $\delta(f)$ is the Dirac delta function.
The right-hand side of \eqref{conv_wave_eqn} represents the source terms for the FWH noise model and are calculated from the velocities and pressure (and its time derivative) that are computed by the compressible Navier-Stokes solver. In this implementation, quadrupoles are neglected and so are their contributions to the viscous stress tensor, but they can be easily added if necessary.

The acoustic pressure $p' = p - p_0$ can be obtained from the density fluctuation, as $p' = c_0^2 \rho '$ in the far field (where $\rho '/ \rho_0 << 1$).
For wind tunnel scenarios, most of the terms needed to solve \eqref{conv_wave_eqn} cancel out, and as a result, only two components of the acoustic pressure fluctuation are solved, namely
\begin{equation} \label{accP_eqn}
    p'(\mathbf{x},t) = p_T'(\mathbf{x},t) + p_L'(\mathbf{x},t),
\end{equation}
\begin{equation} \label{accPT_eqn}
    4\pi p_T' = \int_S{\Bigg[ \frac{\big( 1 - M_{0i}R_i \big) \dot{Q_i}n_i}{R^*}\Bigg]_{\tau_e} dS } \\
    - \int_S{\Bigg[ \frac{U_{0i} \hat{R^*_i} Q_i n_i}{R^{*2}}\Bigg]_{\tau_e} dS },
\end{equation}
\begin{equation} \label{accPL_eqn}
    4\pi p_L' = \int_S{\Bigg[ \frac{\dot{L}_{ij} n_j \hat{R_i}}{c_0 R^*}\Bigg]_{\tau_e} dS } \\
    + \int_S{\Bigg[ \frac{L_{ij}n_j \hat{R^*_i}}{R^{*2}}\Bigg]_{\tau_e} dS },
\end{equation}
where $p_T$ is the thickness pressure fluctuation, $p_L$ is the loading pressure fluctuation, $c_0$ is the speed of sound, $M_0$ is the free-stream Mach number, $U_0$ is the free-stream velocity, $R$ is the phase radius, $R^*$ the amplitude radius (see~\cite{Garrick54atheoretical} for a more detailed explanation), and the hat ( $\hat{}$ ) of both $R$ and $R^*$ represents the partial derivatives of the radius quantities.  The subscript $\tau_e$ indicates that the integrals are calculated at the emission retarded time defined as $\tau_e = t - R/c_0$, see~\cite{Farassat2007DerivationOF} for more details.

For static observers, the retarded time can be calculated analytically, as opposed to when considering bodies with motion. A source time-dominant algorithm is used~\cite{Brentner_2003}, where at each numerical time-step (also known as \emph{emission} time) the \emph{observer} time is calculated, allowing for a time-advancing approach.
All integrals from \eqref{accPT_eqn} and \eqref{accPL_eqn} are numerically evaluated on each h-mesh face of the FWH surface, using the high-order (p-mesh) discretisation at the face, which is inherently related to the Gauss-Legendre (or Gauss-Lobatto) nodes.

%A difficulty doing the reconstruction of the acoustic pressure history using the time-dominant algorithm arises, being related with the fact that for each point on the surface, there are different values of the emission retarded time (equation~\ref{te_eqn}). This creates a non uniformly distribution of the acoustic pressure, i.e. the sound emitted at a certain time (numerical simulation time) will arrive at a different time for each nodal point to the observer. To handle this, a linear interpolation of the face acoustic pressure history is performed (for each receiver) before adding the contribution of each of the faces.

\begin{figure*}[htp]
    \centering
    \caption{Stationary monopole in a moving medium \cite{Lockard_2002}. Comparison between analytical and FWH acoustic pressure.}
        \includegraphics[clip,width=0.8\linewidth]{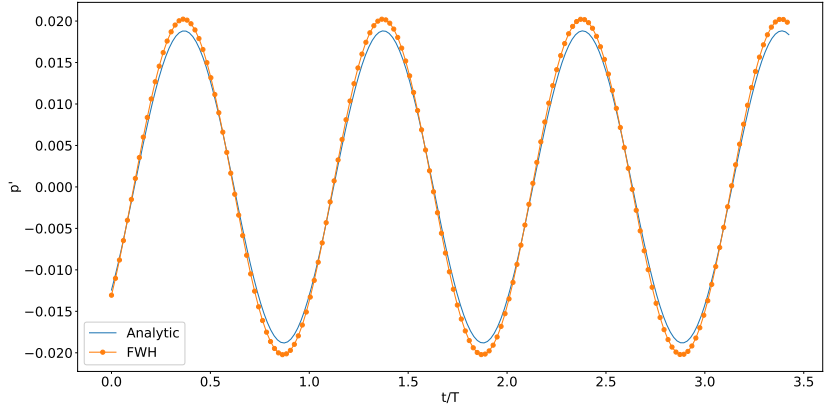}
        \label{fig:fwh_monopole}
\end{figure*}

A classic validation case for evaluating an FWH implementation is the stationary monopole in a moving medium~\cite{Lockard_2002}. The tested monopole has a reference length $l=1$, an amplitude $A/(lc_0)=9.98\times10^{-5}$, an angular frequency $\omega l/c_0=4\pi$, a free-stream Mach number $M_0=0.5$, and a centre at the origin $x_0=[0,0,0]$. Fig.~\ref{fig:fwh_monopole} shows the prediction of the acoustic pressure for an observer located at $500 l$ in the streamwise direction, where a good prediction of the frequency against the analytical solution is obtained, while there is a relative error of around 7\% for the amplitude. 

\begin{figure*}[htp]

    \centering
        \subfloat[Dilation field.]{
        \includegraphics[trim={15cm 0cm 0cm 0cm},clip,width=0.49\linewidth]{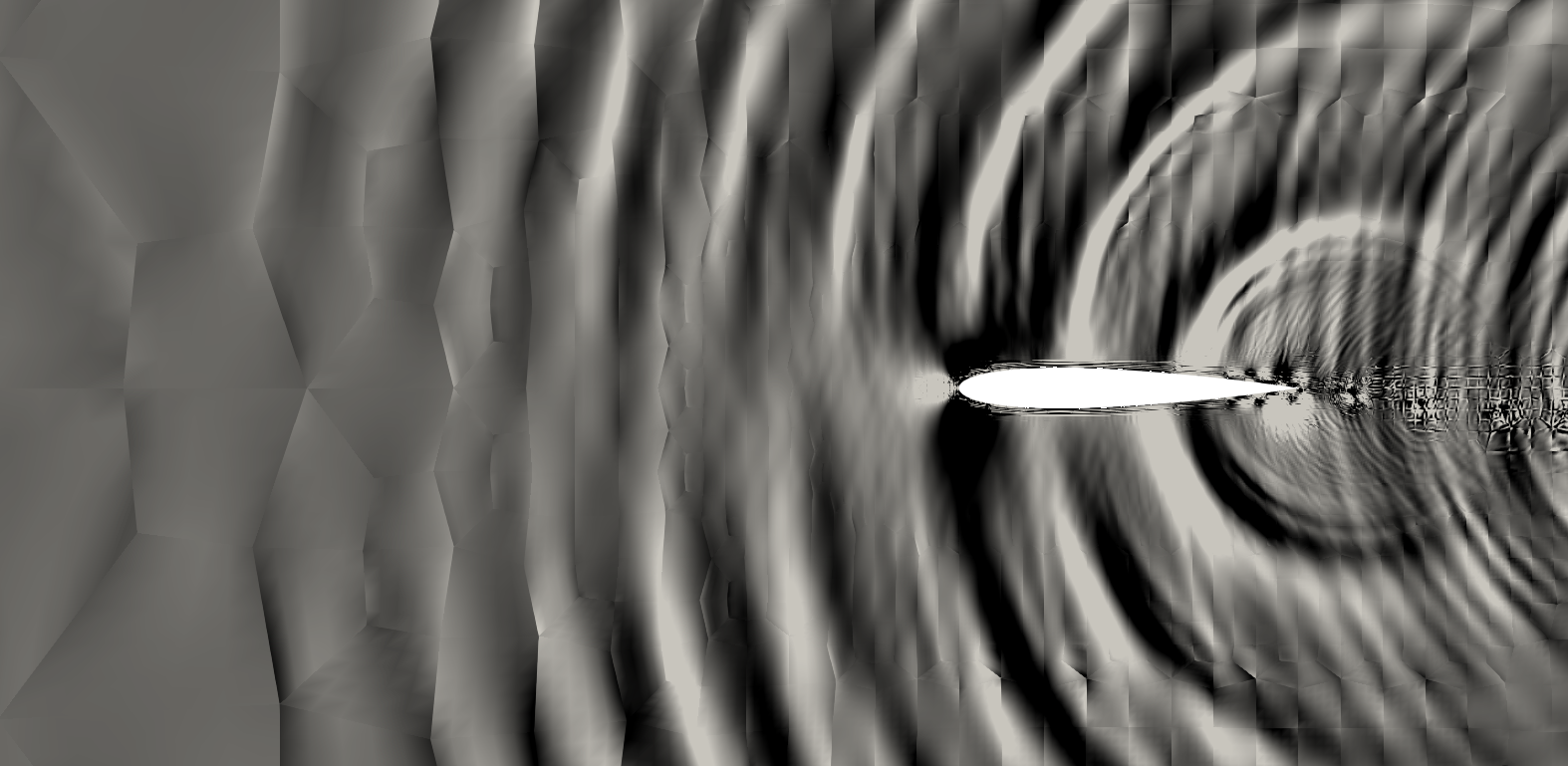}
    \label{fig:Dilation}
    }
    \subfloat[Far field SPL directivity obtained using FWH.]{
        \includegraphics[trim={14,5cm 2cm 13cm 2cm}, clip,width=0.49\linewidth]{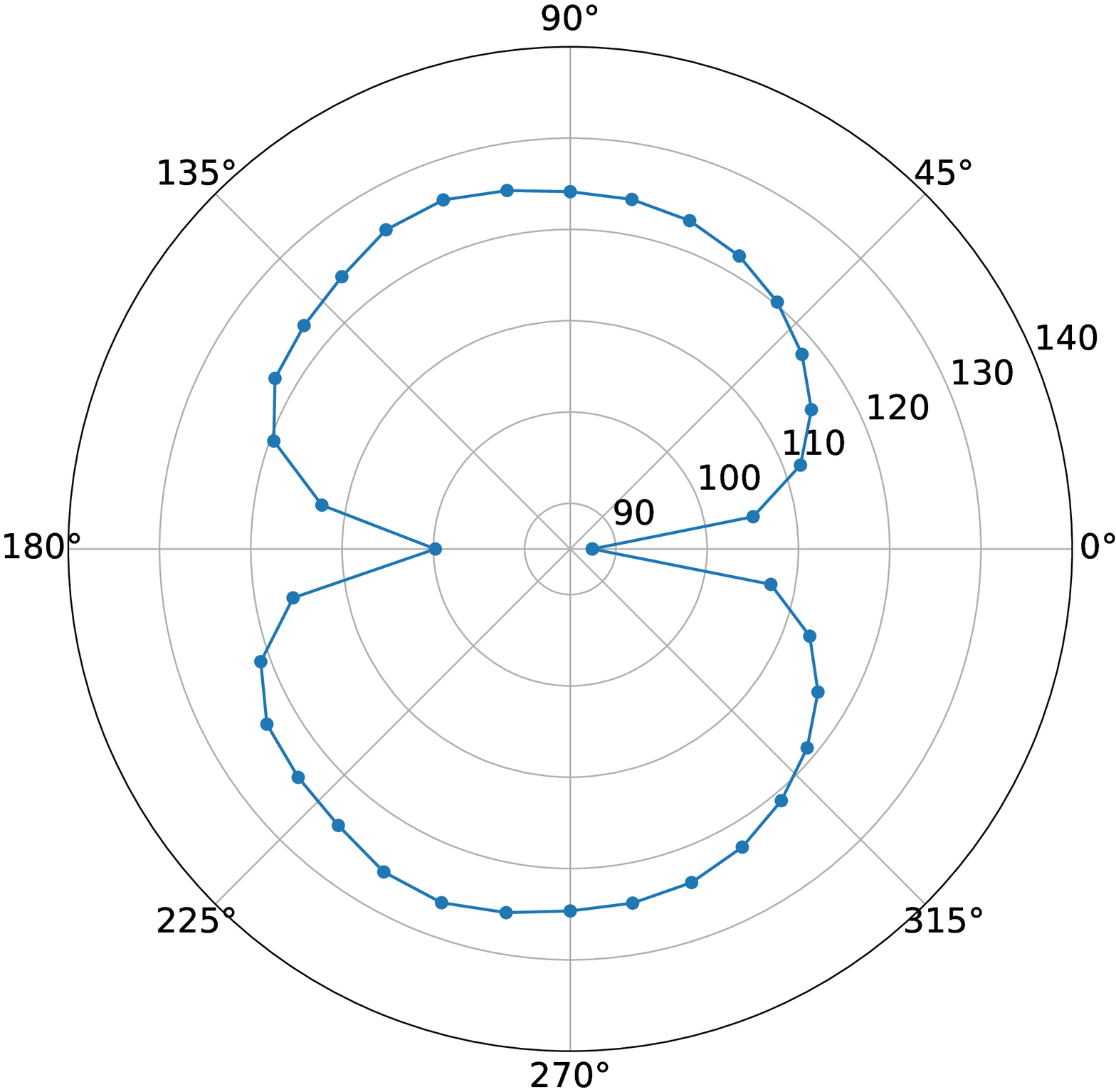}
        \label{fig:fwh_naca}
    }
	\caption{NACA 0012 aeroacoustics at $\text{Re}=1\times10^5$.}
    \label{fig:naca_noise}
\end{figure*}

Another case of classical study is the sound generated by an airfoil in an external flow~\cite{Paterson_1973,DESQUESNES_2007}, where the noise is scattered by the trailing edge. The setup is a NACA 0012 at $\text{Re}=1\times10^5$ based on the airfoil chord, and $\text{Ma}=0.4$. DNC computations can be performed in the near field, as shown in Fig.~\ref{fig:Dilation} for the dilation field, where acoustic waves can be seen to be generated at the trailing edge and propagated to the far field. Furthermore, the directivity of the acoustic pressure in the far field (at a $r=2C$) can be seen in Fig.~\ref{fig:fwh_naca}, computed by the FWH analogy. The expected dipole shape is obtained for this case.

% note from Oscar: I've left the 3 figures, but I can delete one ore more (and the text that comes with it) without problem.

\subsection{Mesh free - Immersed Boundary method}
\label{sec:IB}
In addition to classic body fitted simulations, {\fontfamily{qcr}\selectfont  HORSES3D} includes an Immersed Boundary (IB) technique. IB is a mesh-free method that allows the simulation of complex and moving geometries on simple Cartesian meshes. By avoiding the necessity of body fitting mesh, the grid generation process can be skipped, thus shortening a very time-consuming step.

IB was first introduced by Peskin~\cite{PESKIN1972252} and comes in many flavours; see the review by Kim and Choi \cite{KIM2019301}. %One of the main advantages of this approach is the possibility of simulate arbitrarily complex geometries avoiding the burden of generating a body fitted mesh, which takes a large amount of time and effort. 
A possible implementation of IB is volume penalisation, where to simulate the presence of a geometry, a source term is applied to the Navier-Stokes equations. This approach has recently been extended to high-order methods by our group \cite{KOU2022110817,KOU2022110721}. The method is simple, robust, and can be easily extended to moving geometries.

The IB mask $\chi(\vec{x},t)$ distinguishes the penalised region (inside the geometry) $\Omega_s$ and the outer region (fluid region), $\Omega_f$, as
\begin{equation}
    \chi(\vec{x},t) = 
    \begin{cases}
    1, & \text{if } \vec{x} \in \Omega_s \\
    0, &\text{if } \vec{x} \in \Omega_f,
    \end{cases}
\end{equation}
where $\vec{x}$ is the coordinate vector.
The geometry is considered to be a porous media whose permeability approaches zero to mimic solid geometries. Given the velocity of the geometry $\vec{u}_s = \left( u_s, v_s, w_s \right)$, a source term
\begin{equation}
    \textbf{S} = \frac{\chi}{\eta} 
    \begin{pmatrix}
    0 \\
    \rho u_s - \rho u \\
    \rho v_s - \rho v \\
    \rho w_s - \rho w \\
    \frac{\rho}{2} \left( u_s^2 + v_s^2 + w_s^2 \right) - \frac{\rho}{2} \left( u^2 + v^2 + w^2 \right). \\
    \end{pmatrix}
    \label{S_termIBM}
\end{equation}
 is added to the Navier-Stokes equations, where $\eta$ is the permeability, $\rho$ and $u, v, w$ the density and velocity components, respectively.
Once the mask is generated, the source term \eqref{S_termIBM} is applied to the degrees of freedom and the equations are solved.

In {\fontfamily{qcr}\selectfont  HORSES3D}, the mask is constructed using a ray-tracing approach. The input is a CAD (STL) file representing the geometry to be simulated and a generic mesh. A surface-area heuristic KD-tree is built that contains the whole body. This type of KD-tree is known to improve the efficiency of the ray-tracing algorithm. The construction of the KD-tree is parallelised using MPI and OpenMP. The OpenMP construction implements breadth-first and depth-first parallelisation. The first is used at the starting levels of the KD-tree where the number of triangles, defining the surface geometry in the STL file, is usually big enough to justify the parallelisation of the loops performed on the objects; the latter is used so that each of the remaining KD-tree's branches is performed by a single thread (interested readers are referred to~\cite{wu2011sah}). When MPI is used, the CAD file is split into a number of partitions equal to the number of slaves. A  KD-tree is built for each partition, reducing the computational time since the total number of triangles in the partitions is smaller than the initial one.

A key aspect of IB is the computation of the aerodynamic forces acting on the body, which requires a reconstruction of the surface data. The idea is to find the value of the state variables on the interpolation points over the body surface (IP) and perform surface integration. In general, for what concerns the IB approach, none of the degrees of freedom where the solution is computed lies on the body surface; hence, interpolation is needed. We use the inverse distance weight to compute the state at each interpolation point:
\begin{equation}
    Q_{IP} = \frac{\sum_i\frac{Q_i}{d_i}}{\sum_i\frac{1}{d_i}},
\end{equation}
where $Q_{IP}$ is the state at the interpolation point, $d_i$ is the distance between the $i^{th}$-degree of freedom and the surface normal passing through the interpolation point, and $Q_{i}$ is the state at the $i^{th}$-degree of freedom.
Data reconstruction is performed starting from the so-called band region, i.e. the set of all degrees of freedom lying in the fluid region ($\Omega_f$) close enough to the body surface. These degrees of freedom are candidates for being inside the set of points used for the interpolation. {\fontfamily{qcr}\selectfont  HORSES3D} uses a nearest neighbour algorithm to find the $n$ points (where $n$ is a user-defined parameter) closer to the IP point and, once all the values of the IPs are known, integration is performed. 

\begin{figure*}[htp]
    \centering
            \subfloat[The STL geometry is included for \\
            post-processing clarity.]{
         \includegraphics[width=0.47\textwidth]{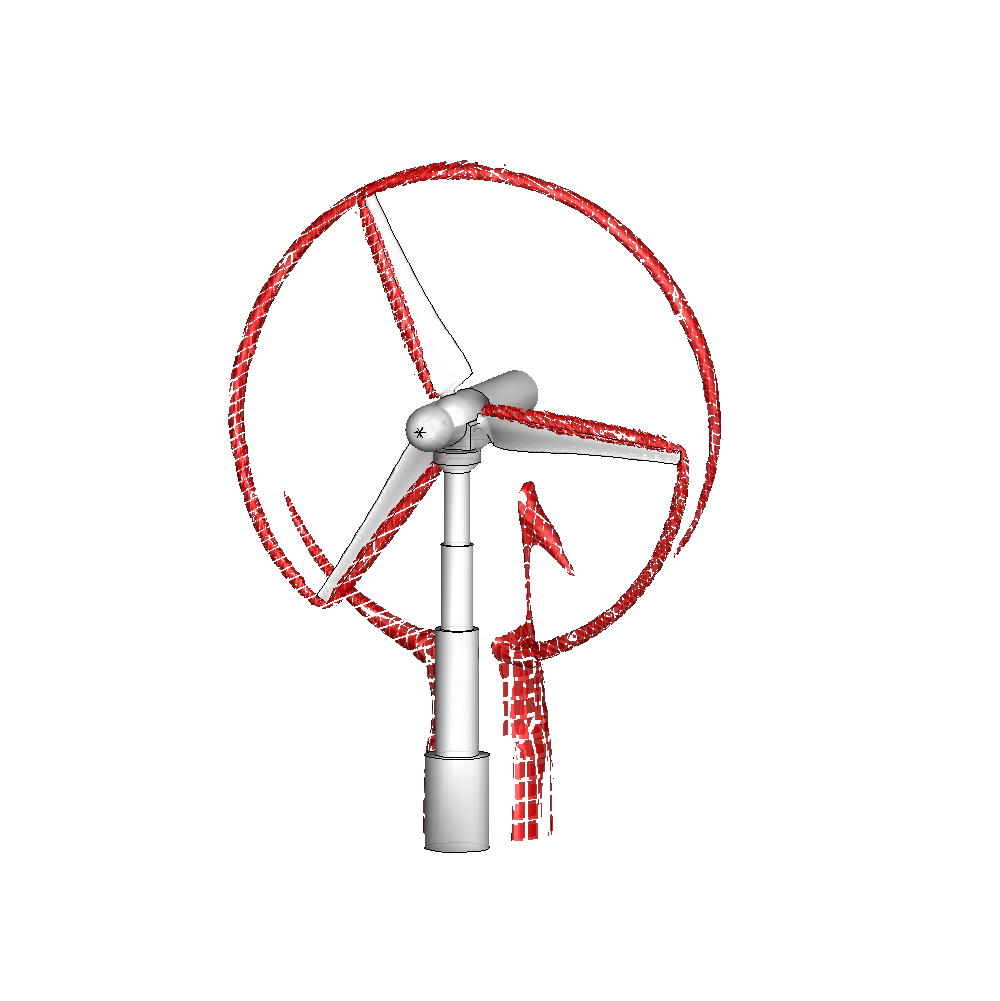}
    \label{fig:VAWT1}
    }
    \subfloat[The high-order mesh is shown. Near the rotor and wake a $P=3$ polynomial is used whilst $P=2$ is selected far from the rotor and wake.]{
         \includegraphics[width=0.47\textwidth]{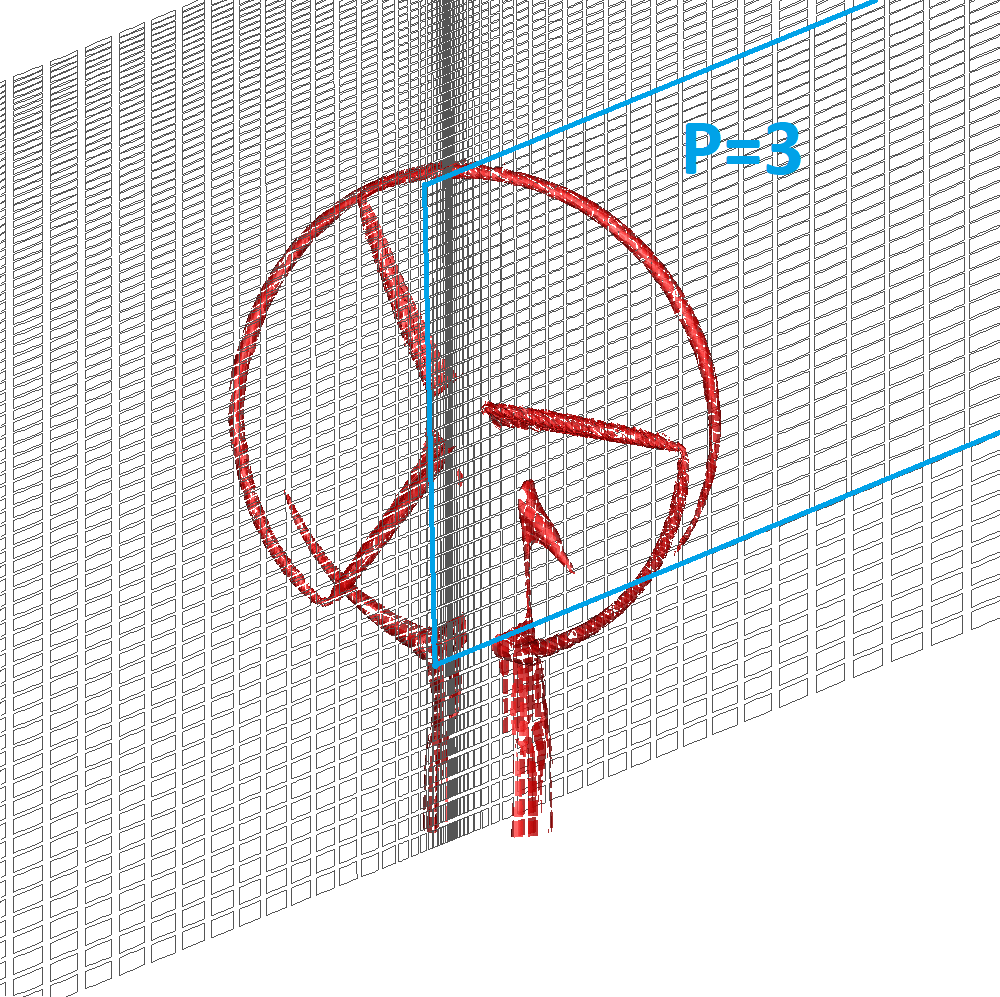}
        \label{fig:VAWT2}
    }
	\caption{Wind turbine simulation with immersed boundary method. Isosurface of velocity magnitude showing tip vortex and the acceleration region at the blade leading edge.} \label{fig:IBM_wt}
\end{figure*}

Fig. \ref{fig:IBM_wt} shows how the IB method can be used to simulate a horizontal wind turbine. In this case the tower and nacelle are static, but the blades rotate with a constant rotational speed. It can be seen that the IB method captures the leading edge accelerations and blade tip vortices without the necessity of generating a complex body-fitted mesh. As an alternative for compute horizontal axis turbines, we have recently implemented various actuator line formulations.

\section{Licensing, code structure, manuals and test cases}
The code is maintained and available on GitHub. 
{\fontfamily{qcr}\selectfont  HORSES3D} is a copyright of NUMATH \hyperlink{numath}{https://numath.dmae.upm.es} under the MIT License. NUMATH (Numerical Methods in Aerospace Technology) is a research group belonging to the School of Aeronautics in Madrid (Universidad Politécnica de Madrid). The objective of the group is to develop high-order numerical methods.

\subsection{GitHub and continuous integration}

The source code of {\fontfamily{qcr}\selectfont  HORSES3D} is stored on GitHub. It has a number of test cases to check the various physics and numerical methods. The objective of test cases is to quickly check if a new commit ``breaks'' the code. To that end, once a new functionality is implemented, a test case is generated that checks it. The solution of this test case is stored and compared to that obtained with new versions of the code. The process of running the test cases and checking the solutions is integrated within the GitHub workflow framework. It is set to run automatically with a pull request trigger. Furthermore, a more complete test case suit is run on a monthly basis to check the code in more detail (e.g., compilation and run with different compilers in serial/parallel, etc.). 

\subsection{Documentation} \label{sec:documentation}

{\fontfamily{qcr}\selectfont  HORSES3D} includes a $doc$ folder with documentation and user manuals where all parameters are identified. The solver uses control files, $*.control$, where the parameters and boundary conditions for the simulations are set.
\begin{tcolorbox}[colback=black!5!white,colframe=gray!75!black]
A variety of test cases are included for the different sets of equations: Cahn-Hilliard, Euler, Incompressible Navier-Stokes, Multiphase, Navier-Stokes, Particles or Navier-Stokes coupled to Cahn-Hilliard. For example, for the compressible Navier-Stokes we include: ActuatorLine, Cylinder, EnergyConservingTest, FlatPlate, ForwardFacingStep, ManufacturedSolutions,  NACA0012, TaylorGreen. 
\end{tcolorbox}
These cases are integrated into the Github workflow framework and run automatically (see previous section). 

\section{Future directions}
We are constantly updating {\fontfamily{qcr}\selectfont  HORSES3D} with new features, multiphysics and coupling it with new techniques. For example, we have recently linked {\fontfamily{qcr}\selectfont  HORSES3D} to neural networks to accelerate simulations and have allowed the solver to call other external libraries through Python interfaces. Future directions include:
\begin{itemize}
    \item Neural networks to accelerate high-order methods:
    
    High-order discontinuous Galerkin methods allow accurate solutions through the use of high-order polynomials inside each mesh element. Increasing the polynomial order leads to high accuracy, but increases the cost. On the one hand, high-order polynomials require more restrictive time-steps when using explicit temporal schemes, and on the other hand, the quadrature rules lead to more costly evaluations per iteration. In a recent work, we propose to accelerate high-order DG methods using Neural Networks.
To this aim, we train a Neural Network using a high-order discretisation, to extract a corrective forcing that can be applied to a low order solution with the aim of recovering high-order accuracy. With this corrective forcing term, we can run a low-order solution (low cost) and correct the solution to obtain high-order accuracy; see details in \cite{DELARA2022105274,DELARAnew3d}.

    \item Coupling with other solvers (e.g., MFEM \cite{mfem-web,mfem}) using Python: 
    
It is often advantageous to link flow solvers with external libraries (e.g., when performing fluid-structure interaction coupling or conjugate heat transfer). In these cases, external libraries can be called to perform part of the simulation (structural calculation or solving heat/radiation). To allow for this flexibility, we have developed an interface from Fortran to Python. Using this interface, it is possible to call Python routines at any point of the CFD calculation. Additionally, it is possible to call other libraries written in C/C++ from Python. 
\end{itemize}

\section{Conclusions}
{\fontfamily{qcr}\selectfont  HORSES3D} is an open-source parallel framework for the solution of non-linear partial differential equations. The solver allows the simulation of compressible flows (with and without shocks), incompressible flows, RANS and LES turbulence models, particle dynamics, multiphase flows, and aeroacoustics.
The numerical schemes provide fast computations, especially when selecting high-order polynomials and MPI/OpenMP combinations for a large number of processors. Entropy-stable schemes, local anisotropic p-adaptation, and efficient dual time-stepping and multigrid allow the solver to tackle problems of industrial relevance, such as aircrafts and wind turbines. 

The modularity and object-oriented programming architecture allow for easy integration of new physics and numerical methodologies, ensuring the necessary  flexibility in a research focused solver. 


\section*{Acknowledgements}
EF, GN and WL acknowledge the financial support of the European Union’s Horizon 2020 research and innovation programme under the Marie Skłodowska-Curie grant agreement (MSCA ITN-EID-GA ASIMIA No 813605). 
GR and EV acknowledge the funding received by the Grant SIMOPAIR (Project No. RTI2018-097075-B-I00) funded by MCIN/AEI/ 10.13039/501100011033 and by ERDF A way of making Europe. 
OM and EV thank the European Union Horizon 2020 Research and Innovation Program under the Marie Sklodowska-Curie grant agreement No 860101 for the zEPHYR project.
SC and EV thank the European Union Horizon 2020 Research and Innovation Program under the Marie Sklodowska-Curie grant agreement No 955923 for the Ssecoid project. DAK is supported by a grant from the Simons Foundation (\#426393, David Kopriva).
Finally, all authors gratefully acknowledge the Universidad Politécnica de Madrid (www.upm.es) for providing computing resources on Magerit Supercomputer.

%In addition, the project has received funding from the European High-Performance Computing Joint Undertaking (JU) under grant agreement (No 956104). The JU receives support from the European Union’s Horizon 2020 research and innovation programme under grant agreement (No 823844) and Spain, France, Germany.

 \bibliographystyle{elsarticle-num} 
 \bibliography{horses3d}





% \end{thebibliography}

\end{document}